\documentclass[11pt,a4paper]{article}

\usepackage{ifthen,latexsym,amssymb,amsmath,bbm,fixmath,stmaryrd}
\usepackage[nobysame,initials]{amsrefs}
\usepackage[shortlabels]{enumitem} % to use: \begin{enumerate}[(a)] \item ... \end{enumerate}
\usepackage{listings} % to format code
\usepackage{hyperref}
\hypersetup{
    colorlinks=true,
    linkcolor=blue,
    filecolor=magenta,      
    urlcolor=cyan,
    }

\usepackage{algorithm2e}
\usepackage{tabularx}
\usepackage{float}
\usepackage{etoolbox}

\newcommand{\e}{\varepsilon}

\newcommand{\C}[1]{{\protect\mathcal{#1}}}
\newcommand{\B}[1]{{\bf #1}}
\newcommand{\I}[1]{{\mathbbm #1}}
\renewcommand{\O}[1]{\overline{#1}}

\newcommand{\V}[1]{\mathbold{#1}}

\newcommand{\hide}[1]{}
\newcommand{\me}{\mathrm{e}}

\newtheorem{definition}{Definition} [section]
\newtheorem{theorem}[definition]{Theorem}
\newtheorem{lemma}[definition]{Lemma}

\newtheorem{corollary}[definition]{Corollary}

\newtheorem{claim}[definition]{Claim}

\usepackage{changepage}

\usepackage{pgf,tikz}
\usepackage{pdflscape}

\newcommand{\bpf}[1][Proof.]{\smallskip\noindent{\it #1} }
\newcommand{\qed}{\nolinebreak\mbox{\hspace{5 true pt}%
  \rule[-0.85 true pt]{3.9 true pt}{8.1 true pt}}}
\newcommand{\epf}{\qed \medskip}

\newcommand{\formatnum}[2]{
  \pgfmathprintnumber[fixed, precision=\ifstrempty{#1}{10}{#1}]{#2}...
}

\begin{document}

%\tableofcontents

\newcommand{\PST}[1]{\ifthenelse{\equal{#1}{}}
{\cite{PikhurkoSliacanTyros19}}
{\cite{PikhurkoSliacanTyros19}*{#1}}}

\newcommand{\LPSS}[1]{\ifthenelse{\equal{#1}{}}
{\cite{LiuPikhurkoSharifzadehStaden23}}
{\cite{LiuPikhurkoSharifzadehStaden23}*{#1}}}

\newcommand{\EL}[1]{\ifthenelse{\equal{#1}{}}
{\cite{EvenzoharLinial15}}
{\cite{EvenzoharLinial15}*{#1}}}

\newcommand{\BGHKS}[1]{\ifthenelse{\equal{#1}{}}
{\cite{BasitGranetHorsleyKundgenStaden25x}}
{\cite{BasitGranetHorsleyKundgenStaden25x}*{#1}}}

\newcommand{\obj}[2][\kappa,\ell]{\lambda_{#1}\ifstrempty{#2}{}{(#2)}}
\newcommand{\Obj}[2][\kappa,\ell]{\Lambda_{#1}\ifstrempty{#2}{}{(#2)}}

\newcommand{\N}{N}
\newcommand{\p}{p}
\renewcommand{\P}{P}
\newcommand{\aut}{\mathrm{aut}}
\newcommand{\tinj}{t}
\renewcommand{\sqcup}{+}
\renewcommand{\mid}{:}

\newcommand{\blow}[2]{\ifstrempty{#2}
{#1{()}}
{#1(#2)}
}

\newcommand{\dedit}{\delta_{\mathrm{edit}}}

\renewcommand{\ge}{\geqslant}
\renewcommand{\geq}{\geqslant}
\renewcommand{\le}{\leqslant}
\renewcommand{\leq}{\leqslant}
\renewcommand{\preceq}{\preccurlyeq}
\renewcommand{\succeq}{\succcurlyeq}

\newcommand{\arxiv}[2]{#1}

%\section{Tables to produce}

\newcommand{\construction}[2]{#1}
\newcommand{\PerfectStability}[2]{#1 &#2}
\newcommand{\MultiPartite}[1]{T_{#1}}

\tikzset{
  vertex/.style={
    circle,
    draw,
    fill=black,
    minimum size=3pt,
    inner sep=0pt
  },
  edge/.style={},
  edgeR/.style={red,line width=1.2pt},
  edgeB/.style={blue,line width=1.2pt}
}

\newcommand{\DrawGraph}[2]{
  \begin{tikzpicture}
    \def\graphradius{10pt} 
    \foreach \i in {0, ..., \the\numexpr #1-1 \relax} {
      \pgfmathsetmacro{\angle}{90 - \i * (360/#1)}
      \node[vertex] (v\i) at (\angle:\graphradius) {};
    }
    \foreach \u / \v in {#2} {
      \draw[edge] (v\u) -- (v\v);
    }
  \end{tikzpicture}
}

\usetikzlibrary{calc}

\newcommand{\DrawSemiGraph}[3]{
  \begin{tikzpicture}[baseline={($(current bounding box.center)-(0pt,3pt)$)}]
    \def\graphradius{7pt} 
    \foreach \i in {0, ..., \the\numexpr #1-1 \relax} {
      \pgfmathsetmacro{\angle}{90 - \i * (360/#1)}
      \node[vertex] (v\i) at (\angle:\graphradius) {};
    }
    \foreach \u / \v in {#2} {
      \draw[edgeR] (v\u) -- (v\v);
    }
    \foreach \u / \v in {#3} {
      \draw[edgeB] (v\u) -- (v\v);
    }
  \end{tikzpicture} % 
}

\title{Some exact inducibility-type results for graphs via flag algebras}

\author{Levente Bodn\'ar\thanks{Email: bodnalev@gmail.com%%levente.bodnar@warwick.ac.uk
}\ \mbox{ and }Oleg Pikhurko\thanks{Email: O.Pikhurko@warwick.ac.uk}\\
Mathematics Institute and DIMAP\\
 %Centre for Discrete Mathematics and its Applications (DIMAP)\\
University of Warwick\\
Coventry CV4 7AL, UK}

\maketitle

\begin{abstract} 
The \emph{$(\kappa,\ell)$-edge-inducibility problem} asks for the maximum number of $\kappa$-subsets inducing exactly $\ell$ edges that a graph of given order~$n$ can have. Using flag algebras and stability approach, we resolve this problem for all sufficiently large $n$ (including a description of all extremal and almost extremal graphs) in eleven new non-trivial cases when $\kappa\le 7$. 

We also compute the \emph{$F$-inducibility constant} (the asymptotically maximum density of induced copies of $F$ in a graph of given order $n$) and obtain some corresponding structure results for three new graphs $F$ with $5$ vertices: the 3-edge star plus an isolated vertex, the 4-cycle plus an isolated vertex, and the 4-cycle with a pendant edge.
\end{abstract}

\section{Introduction}

Let $\kappa$ and $\ell$ be given non-negative integers  with $\ell\le \binom{\kappa}2$. A \emph{$(\kappa,\ell)$-graph} is a graph having exactly $\kappa$ vertices and $\ell$ edges. For a graph $G$, let  $\Obj{G}$ denote the number of \emph{$(\kappa,\ell)$-subgraphs}, meaning induced subgraphs of $G$ with $\kappa$ vertices and $\ell$ edges. In other words, $\Obj{G}$ counts $\kappa$-subsets of $V(G)$ that span exactly $\ell$ edges in~$G$. The \emph{$(\kappa,\ell)$-edge-inducibility problem} (or \emph{$(\kappa,\ell)$-problem} for short) asks for 
\[
\Obj{n}:=\max\{\Obj{G}\mid \mbox{$n$-vertex graph $G$}\},
\]
 the maximum number of $(\kappa,\ell)$-subgraphs that a graph with $n$ vertices can have. It is natural to consider the normalised function
 \[
 \obj{n}:=\frac{\Obj{n}}{\binom{n}{\kappa}},\quad \mbox{for $n\ge \kappa$}.
 \]
 The standard averaging argument shows that the limit
 \begin{equation}\label{eq:lim}
 \obj{}:=\lim_{n\to\infty} \obj{n}
 \end{equation}
 exists, see e.g.~\PST{Lemma~2.2}. 
 %Following \cite{AlonHefetzKrivelevichTyomkyn20} 
 We refer to the value of the limit as  
 %\emph{edge-inducibility} (or the \emph{edge-inducibility constant})
 the \emph{edge-inducibility constant} 
 of~$(\kappa,\ell)$.

 Observe that, by replacing all graphs with their complements, the value of $\obj{n}$ will not change if we replace $\ell$ by $\binom{\kappa}{2}-\ell$.
 %\Obj{n}=\Obj[\kappa,\binom{\kappa}{2}-\ell]{n}$.
 %(and thus e.g.\
 % $\obj{n}=\obj[\kappa,\binom{\kappa}{2}-\ell]{n}$ and 
 %$\obj{}=\obj[\kappa,\binom{\kappa}{2}-\ell]$). 
 Trivially, it holds that $\obj[\kappa,\ell]{n}=1$ if and only if $\ell=0$ or $\binom{\kappa}2$. Also, the special case $\ell=1$, which is equivalent to the inducibility problem  for the unique up to isomorphism $(\kappa,1)$-graph, is resolved through the results in~\cite{BrownSidorenko94,Hirst14,LiuMubayiReiher23,LiuPikhurkoSharifzadehStaden23,LiuMaZhu}, where in particular Liu, Mubayi and Reiher~\cite[Theorem 1.13]{LiuMubayiReiher23} obtained an explicit formula for $\obj[k,1]{}$ valid for every~$\kappa\ge 4$.
 Thus we restrict ourselves to $2\le \ell\le \binom\kappa2/2$ only.

This problem was recently introduced by Alon, Hefetz, Krivelevich and Tyomkyn~\cite{AlonHefetzKrivelevichTyomkyn20} and has received considerable attention. In particular, the Edge-Statistics Conjecture~\cite[Conjecture~1.1]{AlonHefetzKrivelevichTyomkyn20} that %(excluding the trivial cases when $\ell$ is $0$ or $\binom\kappa2$) 
$\obj{}\le 1/\me+o_{\kappa}(1)$ was fully resolved by a sequence of papers by Kwan, Sudakov and Tran~\cite{KwanSudakovTran19}, Martinsson,
Mousset, Noever and Truji\'c~\cite{MMNT19} and Fox and Sauerman~\cite{FoxSauerman20}. Two other conjectures of Alon et al~\cite[Conjectures~6.1 and~6.2]{AlonHefetzKrivelevichTyomkyn20} on stronger bounds in the case when $\ell$ is well separated from $0$ 
%or $\binom{\kappa}2$, 
were resolved by Kwan, Sudakov and Tran~\cite{KwanSudakovTran19} and Kwan and Sauermann~\cite{KwanSauermann23x}.  Hypergraph versions of these results were very recently obtained by Jain, Kwan, Mubayi and Tran~\cite{JainKwanMubayiTran25x}.  Motivated by these results, versions of this problem for hypercubes were studied by Alon, Axenovich and Goldwasser~\cite{AlonAxenovichGoldwasser24} and the authors~\cite{BodnarPikhurko25}.

Here, we systematically investigate the cases when $\kappa\le 7$ of the edge-inducibility problem for graphs using the flag algebra method. Our aim is not only to find the value of the edge-inducibility constant $\obj{}$ (that is, to determine $\Obj{n}$ within an additive $o(n^{\kappa})$ error term as $n\to\infty$) but also to prove perfect stability. The formal  definition of this property will appear in Section~\ref{se:PerfectStab}; informally speaking, here it means that there is a blowup pattern $B$ and a constant $C$ such that every graph $G$ with $n\ge C$ vertices can be made into a blowup of $B$ by changing at most $C(\obj{n}-\obj{G})n^{2}$ adjacencies. In particular, every $n$-vertex graph $G$ which maximizes the number of $(\kappa,\ell)$-subgraphs is a blowup of $B$; thus the determination of $\obj{n}$ (and of the set of extremal graphs) amounts to finding optimal part sizes of a $B$-blowup (that is, to maximising some polynomial over non-negative integers summing up to $n$). If this analytic problem is resolved with a description of all optimal part ratios, then perfect stability implies Erd\H os--Simonovits stability~\cite{Erdos67,Simonovits68} 
that aims to describe the structure of every graph $G$ of order $n\to\infty$ with $\obj{G}=\obj{}+o(1)$ up to $o(n^2)$ adjacencies (see Section~\ref{se:PerfectStab}). %for definition and discussion).

We were able to determine the value of $\obj{}$ for eleven new pairs $(\kappa,\ell)$, also showing that perfect stability holds in each solved case except for $(\kappa,\ell)=(4,3)$. Table~\ref{ta:known} summarises our new findings, where we use the following notation for constructions: $\MultiPartite{n_0,\dots,n_{m-1}}$ is the complete $m$-partite graph with parts of sizes $n_0,\dots, n_{m-1}$; $K_n$ is the clique with $n$ vertices, $F\sqcup H$ is the vertex-disjoint union of  graphs $F$ and $H$, $mF:=F\sqcup\dots\sqcup F$ is the union of $m$ copies of $F$, and the constants are 
\begin{align}
    \alpha_1 &:= \frac{1}{16} \left(9 - \sqrt{17} \right),\label{eq:a1} \\
    \alpha_2 &:= \frac{1}{2} \left(1-\sqrt{\frac{1}{3} \left(2 \sqrt{10}-5\right)}\right), \label{eq:a2}\\
    \alpha_3 &:= \frac{1}{2} \left(1 - \sqrt{\frac13 \left( \frac45 \sqrt{10} - 1 \right)} \right).\label{eq:a3}
\end{align}
We refer the reader to Section~\ref{se:computer} for the formal statements and further details.

% it is rather squeezed in normally, this adjusts the separation
\renewcommand{\arraystretch}{1.5}

%\lb{not sure what is the standard notation for complete graph blowup is, https://arxiv.org/pdf/1312.1205 uses $K_{a, b, c}$, but I think it is not that consistent with the cliques, so I use T with a macro for now}

\begin{table}[h!]
\begin{center}
\begin{tabular}{c|c|c|c|c}
\hline
$(\kappa,\ell)$ & Construction & $\obj[\kappa,\ell]{}$ & Stability & Reference\\
\hline
\hline 

%$(3,1)$ & \construction{$2K_{\frac{n}{2}}$}{\cite{BrownSidorenko94}} & $\frac34$ & \PerfectStability{Perfect stability}{follows from \cites{BrownSidorenko94,LiuPikhurkoSharifzadehStaden23}} \\
%\hline

%$(4,1)$ & \construction{$5K_{\frac{n}{5}}$}{\cite{BrownSidorenko94}} & $\frac{72}{5^3}$ & \PerfectStability{Perfect stability}{follows from \cites{BrownSidorenko94,LiuPikhurkoSharifzadehStaden23}} \\

$(4,2)$ & \construction{%$2\MultiPartite{\frac{n}{6}, \frac{n}{6}, \frac{n}{6}}$
$2\MultiPartite{{n}/{6},\, {n}/{6},\, {n}/{6}}$}{\ref{th:constrs}} & $\frac12$ & \PerfectStability{Perfect stability}{Theorem~\ref{th:stats_42}}\\

$(4,3)$ & \construction{%$2K_{\frac{n}{2}}$ or $\MultiPartite{\frac{n}{2}, \frac{n}{2}}$
$2K_{{n}/{2}}$ or $\MultiPartite{{n}/{2},\, {n}/{2}}$}{\ref{th:constrs}} & $\frac12$ & \PerfectStability{See Section~\ref{se:43}}{Theorem~\ref{th:stats_43}}\\

\hline

%$(5,1)$ & \construction{$8K_{\frac{n}{8}}$}{\cite{BrownSidorenko94}} & $\frac{525}{2^{10}}$ & \PerfectStability{Unknown}{}\\

$(5,2)$ & \construction{%$3\MultiPartite{\frac{n}{9}, \frac{n}{9}, \frac{n}{9}}$
$3\MultiPartite{{n}/{9},\, {n}/{9},\, {n}/{9}}$}{\ref{th:constrs}} & $\frac{280}{3^6}$ & \PerfectStability{Perfect stability}{Theorem~\ref{th:stats_52}}\\

$(5,3)$ & \construction{$\MultiPartite{\alpha_1 n,\, \alpha_1 n} \sqcup K_{(1-2\alpha_1) n}$}{\ref{th:constrs}} & $\frac{255 \sqrt{17} - 535}{2^{10}}$ & \PerfectStability{Perfect stability}{Theorem~\ref{th:stats_53}}\\

$(5,4)$ & \construction{%$2K_{\frac{n}{2}}$
$2K_{{n}/{2}}$}{\ref{th:constrs}} & $\frac58$ & \PerfectStability{Perfect stability}{Theorem~\ref{th:stats_54}}\\

\hline

%$(6,1)$ & \construction{$13K_{\frac{n}{13}}$}{\cite{BrownSidorenko94}} & $\frac{178200}{13^5}$ & \PerfectStability{Unknown}{}\\

$(6,4)$ & \construction{%$3K_{\frac{n}{3}}$
$3K_{{n}/{3}}$}{\ref{th:constrs}} & $\frac{40}{3^4}$ & \PerfectStability{Perfect stability}{Theorem~\ref{th:stats_64}}\\

$(6,5)$ & \construction{$\MultiPartite{\alpha_2 n,\, (1-\alpha_2) n}$}{\ref{th:constrs}} & $\frac{10 \sqrt{10}-28}{9}$ & \PerfectStability{Perfect stability}{Theorem~\ref{th:stats_65}}\\

$(6,7)$ & \construction{%$2K_{\frac{n}{2}}$
$2K_{{n}/{2}}$}{\ref{th:constrs}} & $\frac{15}{2^5}$ & \PerfectStability{Perfect stability}{Theorem~\ref{th:stats_67}}\\

\hline

%$(7,1)$ & \construction{$19K_{\frac{n}{19}}$}{\cite{BrownSidorenko94}} & $\frac{21591360}{19^6}$ & \PerfectStability{Unknown}{}\\

$(7,6)$ & \construction{$\MultiPartite{\alpha_3 n,\, (1-\alpha_3) n}$}{\ref{th:constrs}} & $\frac{28\sqrt{10} - 35}{135}$ & \PerfectStability{Perfect stability}{Theorem~\ref{th:stats_76}}\\

$(7,9)$ & \construction{%$2K_{\frac{n}{2}}$
$2K_{{n}/{2}}$}{\ref{th:constrs}} & $\frac{35}{2^6}$ & \PerfectStability{Perfect stability}{Theorem~\ref{th:stats_79}}\\

$(7,10)$ & \construction{%$\MultiPartite{\frac{n}3, \frac{2n}3}$
$\MultiPartite{{n}/3,\, {2n}/3}$}{\ref{th:constrs}} & $\frac{28}{3^4}$ & \PerfectStability{Perfect stability}{Theorem~\ref{th:stats_710}}\\
\hline
\end{tabular}
\caption{New values of edge-inducibility constants.}
\label{ta:known}
\end{center}
\end{table}

The case $(\kappa,\ell)=(4,3)$ turned out to be special in many aspects (see Section~\ref{se:43} for details). This problem is self-complementary and there are two types of extremal graphs, namely, complete bipartite graphs or two disjoint cliques (with part sizes $n/2+o(n)$). 
%The optimal part sizes are $\frac n2+o(n)$ and 
Interestingly, when we change the adjacency of a pair $xy$ between the parts then the number of $(4,3)$-subgraphs changes by only $O(n)$, instead of a positive fraction of ${n-2\choose 2}=\Theta(n^2)$ $4$-sets containing $xy$. (In fact, if we start with parts of sizes exactly $\lfloor n/2\rfloor$ and $\lceil n/2\rceil$ then changing one adjacency across  strictly increases the number of $(4,3)$-subgraphs.) 
%This can be amplified by changing all pairs between two sets of size $o(n)$ inside each part. 
It follows that perfect stability does not hold for this problem, even if we consider the weaker version where we allow finitely many possible patterns. However, we can prove a version of Erd\H os--Simonovits stability in Lemma~\ref{lm:43stab}: if an $n$-vertex graph $G$ satisfies $\obj[4,3]{G}=\obj[4,3]{}+o(1)$ then $G$ is $o(n^2)$-close to one of the above two constructions. The optimal part sizes can be computed exactly and they deviate from $n/2$ by $\sqrt{3n}/2+O(1)$. This extra imbalance adds a ``drift" that penalises wrong pairs across and suffices for us to prove in Theorem~\ref{th:43exact} that every extremal graph of large order $n$ has no wrong pairs at all.

We also obtained new results on the following graph inducibility problem. (See also Section~\ref{se:concluding} for a discussion on the related semi-inducibility problem.)
For graphs $F$ and $G$ with $\kappa\le n$ vertices respectively, let $\P(F,G)$ denote the number of $\kappa$-subsets of $V(G)$ that span a graph isomorphic to $F$ and let $\p(F,G):=\P(F,G)/\binom n{\kappa}$ be the \emph{density} of $F$ in~$G$.  The \emph{inducibility
problem} for a graph $F$ asks for $\obj[F]{n}$,  the maximum of $\obj[F]{G}:=\p(F,G)$ over all graphs $G$ with $n$ vertices. As before, it is easy to show that the limit $\obj[F]{}:=\lim_{n\to\infty} \obj[F]{n}$ exists; we call it the \emph{inducibility constant} of~$F$.

The inducibility problem has drawn a great amount of interest since it was introduced by Pippenger and Golumbic~\cite{PippengerGolumbic75} in 1975. For some sample of results, see e.g.\ 
\cite{Siran84,BollobasNaraTachibana86,BrownSidorenko94,BollobasEgawaHarrisJin95,Hirst14,HatamiHirstNorin14jctb,%
EvenzoharLinial15,BaloghHuLidickyPfender16,HefetzTyomkyn18,KralNorinVolec19,PikhurkoSliacanTyros19,Yuster19,%
LidickyMattesPfender23,Ueltzen24x}.
% Also, it was considered for many other types of structures:  oriented graphs~\cite{Huang14,ChoiLidickiPfender20,BozykGrzesikKielak22,HuMaNorinWu24}, tournaments~\cite{BurkeLidickyPfenderPhillips21x}, hypercubes~\cite{GoldwasserHansen21,GoldwasserHansen24,BodnarPikhurko25}, binary and $d$-ary trees~\cite{CzabarkaDossouSzekelyWagner17,CzabarkaDossouSzekelyWagner20}, tournaments~\cite{BurkeLidickyPfenderPhillips21x}, 
%constrained graphs (eg max $C_5$ in triangle-free, Erdos $f(k,l)$-problem), in planar graph, etc.

If $F$ is complete partite then the result by Brown and Sidorenko~\cite{BrownSidorenko94}
%~\cite{SchelpThomason98} 
implies that, in order to determine the value of $\obj[F]{n}$, it is enough to consider complete partite graphs on $[n]$ and the problem in the limit reduces to some analytic-type optimisation on the space of part ratios. If the latter is fully solved, with the description of all extremal ratios, then the method of Liu, the second author, Sharifzadeh and Staden~\cite{LiuPikhurkoSharifzadehStaden23} can often be applied to decide if perfect stability holds or not.
So we exclude complete partite $F$ from our consideration. Since the inducibility constant $\obj[F]{}$ does not change if we replace $F$ by its complement, it is enough to consider only one graph from each complementary pair.
% (In particular, the inducibility problem for a graph $F$ which is a union of cliques reduces to 

Each 3-vertex graph $F$ or its complement is complete partite, so we exclude these (as they are covered by the above result of Brown and Sidorenko~\cite{BrownSidorenko94}).

All 4-vertex graphs $F$ were resolved by the results in~\cite{BollobasNaraTachibana86,Exoo86,BrownSidorenko94,Hirst14} except when $F=P_4$ is the 4-vertex path. The best known lower bound $\obj[P_4]{}\ge 1173/5824 = 0.2014...$ is due to~Even-Zohar and Linial~\EL{} while the best known upper bound $\obj[P_4]{}\le 0.204513...$ comes from flag algebras. 

Here we look at $5$-vertex graphs~$F$. For notational convenience, we assume by default that the vertex set of $F$ is $\{0,1,2,3,4\}$ while $xy$ means a pair $\{x,y\}$.
%, whose vertex set by default is $\{1,2,3,4,5\}$.
Even-Zohar and Linial~\EL{Table~2} produced a summary of known and new results for $5$-vertex graphs $F$, in particular providing numerical upper bounds on $\obj[F]{}$ coming from flag algebra calculations in the cases when the exact value was not known.

After the appearance of~\EL{}, two new $5$-vertex cases (when $F$ is not complete partite), namely when $E(F)=\{01,12\}$  
(the 2-edge path plus 2 isolated vertices) and 
$E(F)=\{01,12,23,24\}$ (the ``$Y$-graph'') 
%$E(F)=\{\{1,2\},\,\{2,3\},\, \{3,4\},\, \{3,5\}\}$ (the ``$Y$-graph'') and $E(F)=\{\{1,2\},\,\{2,3\}\}$ 
 were fully resolved for all large $n$ (including perfect stability) by the second author, Sliacan and Tyros~\cite{PikhurkoSliacanTyros19}. %Here we assume that $V(F)=\{a,b,c,d,e\}$ and abbr
The authors of~\cite{PikhurkoSliacanTyros19} also tried to solve some other open $5$-vertex cases using Emil Vaughan's package \texttt{flagmatic} but were not able to. 

Using the first author's new package \texttt{FlagAlgebraToolbox}, see~\cite{Bodnar26fat}, we are able to determine the value of $\obj[F]{}$ for three new 5-vertex graphs $F$: 
%$\MultiPartite{3,1}\sqcup K_1$, 
the 3-star $\MultiPartite{3,1}$ plus an isolated vertex,
the 4-cycle $\MultiPartite{2,2}$ plus an isolated vertex, and the 4-cycle $\MultiPartite{2,2}$ with a pendant edge attached. These results are summarised in Table~\ref{ta:NewF}, where 
%$R_{p}$ denotes a sequence of $p$-quasirandom graphs with orders tending to infinity and 
we let $\beta:=(3+\sqrt{3})/12$ while $R(G,p)$ means a ``typical'' 
spanning 
subgraph of $G$ when each edge of $G$ is kept with probability~$p$ independently of the other edges. (See Section~\ref{se:computer} for the formal statements and further details.)
%\[\beta:=\frac{3+\sqrt{3}}{12}.\]

%We refer the reader to Section~\ref{se:computer} for further details.
\begin{table}[h!]
\begin{center}
\begin{tabular}{c|c|c|c|c}
\hline
Edges of $F$ & Construction & $\obj[F]{}$ & Stability & Reference\\
\hline
\hline 
$01,02,03$ & $\MultiPartite{\beta n,\,\beta n}\sqcup \MultiPartite{(1/2-\beta) n,\,(1/2-\beta)n}$ & $\frac5{24}$ & Perfect stability & Theorem~\ref{th:ind0}\\
$01,12,23,30$ & $\MultiPartite{\beta n,\,\beta n}\sqcup \MultiPartite{(1/2-\beta) n,\,(1/2-\beta)n}$ & $\frac5{32}$ & Perfect stability & Theorem~\ref{th:4Cycle+v}\\
$01,12,23,30,04$ & 
%bipartite ${5/6}$-quasirandom 
$R(\MultiPartite{n/2,\,n/2},\frac 56)$
& $\frac{5^6}{2^8\cdot 3^5}$ & 
%bipartite $\frac56$-quasirandom 
%Within $o(n^2)$ edits
See Section~\ref{se:LeafedC4}
%Theorems~\ref{th:4Cycle+eStab} and~\ref{th:C4+eExtr}
& Theorem~\ref{th:4Cycle+e}\\
\hline
\end{tabular}
\caption{New values of the inducibility constant $\obj[F]{}$ for a graph $F$ with $V(F)=\{0,1,2,3,4\}$.}
\label{ta:NewF}
\end{center}
\end{table}

Note that, in the first two cases, the structure of large extremal graphs (which are disjoint unions of 2 complete bipartite graphs) and the limiting part ratios  happen to be the same. Also, we can prove perfect stability in both cases. 

In the last case when $F$ is the 4-cycle with a pendant edge, we prove in Theorem~\ref{th:4Cycle+eStab} via extra arguments that every almost extremal graph $G$ of order $n\to\infty$ admits a balanced vertex partition $V(F)=V_0\cup V_1$ such that each part spans $o(n^2)$ edges while the induced bipartite graph $G[V_0,V_1]$ is $5/6$-quasirandom, thus obtaining a good characterisation of almost extremal graphs. This is in an interesting contrast with the result of Jain, Michelen and Wei~\cite{JainMichelenWei} that (non-bipartite) Erd\H os--R\'enyi random graphs $R(K_n,p)$ for constant $p\in (0,1)$ cannot be almost extremal for any graph inducibility problem.
Regarding extremal graphs, that is, $n$-vertex graphs $G$ with $\obj[F]{G}=\obj[F]{n}$,  we additionally prove in
Theorem~\ref{th:C4+eExtr} that, for all large $n$, each such graph $G$ admits a vertex partition $V(G)=V_0\cup V_1$ into two independent sets and every vertex of $G$ has $(5/12+o(1))n$ neighbours in the other part.
%(in addition to the bipartite graph $G[V_0,V_1]$ being close to $5/6$-quasirandom). 
This reduces the $F$-inducibility problem for large $n$ to its bipartite version (modulo the issue of finding optimal part sizes). Resolving this bipartite problem exactly seems challenging and we limit ourselves to the above partial description of extremal graphs.

\section{Preliminaries}\label{se:Prelim}

In this section, we present here some definitions and auxiliary results.

Let $\I R$ denote the set of reals. Let $\I N$ denote the set of non-negative integers and, for $n\in\I N$, we define $[n]:=\{0,\dots,n-1\}$. Note that we start indexing from $0$, merely to be consistent with the same convention as in our code. 
If the meaning is clear, we may abbreviate an unordered pair $\{u,w\}$ as $uw$, including the case when $u$ and $w$ are single-digit numbers. For a set $X$ and an integer $\kappa\ge 0$, the set of all $\kappa$-subsets of $X$ is denoted by $\binom{X}{\kappa}$. 
%A pair $\{u,w\}$ may be abbreviated to $uw$. 
Also, $a=b\pm \e$ for $a,b\in\I R$ means $b-\e \le a\le b+\e$. 
%Let $n^{(\kappa)}:=\prod_{i=1}^{\kappa-1}(n-i)$ denote the $\kappa$-th falling power of~$n$. 
We may omit ceiling/floor signs when they are not essential.

%For a graph $G=(V(G),E(G))$, its \emph{order} is $v(G):=|V(G)|$.
% and the \emph{degree} $\deg_G(x)$ of a vertex $x$ is the number of edges incident to $x$. 

A \emph{pattern} $B=(V(B),E(B))$ is a graph where we additionally allow loops on some vertices (but we do not allow multiple edges). Its \emph{order} is $v(B):=|V(B)|$. We write $\{u,u\}\in E(B)$ (or $uu\in E(B)$) to indicate that there is a loop on a vertex $u$. A \emph{pattern automorphism} is a bijection $f:V(B)\to V(B)$ such that, for every $u,w\in V(B)$, we have $\{u,w\}\in E(B)$ if and only if $\{f(u),f(w)\}\in E(B)$; thus it is an automorphism of the underlying graph that also preserves loops and non-loops.

The \emph{$B$-neighbourhood} $\Gamma_B(u)$ of a vertex $u\in V(B)$ is the set $\{w\in V(B)\mid \{u,w\}\in E(B)\}$. Note that $u$ itself is included into $\Gamma_B(u)$ if and only if $u$ is a loop in $B$. The \emph{degree} of $u$ is $\deg_B(u):=|\Gamma_B(u)|$. Of course, these definitions also apply to graphs (which are patterns without loops). Also, for a graph $F$, its \emph{complement} is
$\O F:=\left(V(F),\binom{V(F)}{2}\setminus E(F)\right)$. For $X\subseteq V(F)$, the subgraph \emph{induced} by $X$ is $F[X]:=(X,\{uw\in E(F)\mid u,w\in X\})$ and, for disjoint $X,Y\subseteq V(F)$, we denote $F[X,Y]:=\{(u,w)\in X\times Y\mid uw\in E(F)\}$.

Recall that $K_n$ denotes the complete graph with $n$ vertices, $P_n$ is the $n$-vertex path, and $\MultiPartite{n_0,\dots,n_{m-1}}$ denotes the complete $m$-partite graph with parts of sizes $n_0,\dots,n_{m-1}$. We may refer to $\MultiPartite{2,1}\cong P_3$ as the \emph{cherry}. Also, $F\sqcup H$ denotes the union of vertex-disjoint copies of graphs $F$ and~$H$.
When we define a small graph/pattern, we may write it as $(m,E)$, meaning that the vertex set is $[m]$.
For example, we denote the 4-vertex path  as $(4,\{01,12,23\})$ and 2 isolated loops as $(2,\{00,11\})$. 
%Usually, we try to use the same labelling as in the package (which is the one returned by \textt{nauty}). 

Let $B$ be a pattern with vertex set $[m]$. 
For pairwise disjoint sets $V_0,\dots,V_{m-1}$ (with some possibly empty), the \emph{blowup} $\blow{B}{V_0,\dots,V_{m-1}}$ of $B$ is the graph (without loops) on $V=\cup_{i=0}^{m-1} V_i$ where distinct $x\in V_i$ and $y\in V_j$ are adjacent if and only if $\{i,j\}\in E(B)$. In particular, a part $V_i$ spans a clique (resp.\ an independent set) if $i$ is  (resp.\ is not) a loop of~$B$. Let $\blow{B}{}$ denote the family of all blowups of~$B$. 
A \emph{homomorphism} from a graph $F$ to a pattern $B$
% (possibly with loops) 
is a (not necessarily injective) function $f:V(F)\to V(B)$ such that for every distinct $x,y\in V(F)$ it holds that $\{x,y\}\in E(F)$ if and only if $\{f(x),f(y)\}\in E(B)$. Thus homomorphisms from $F$ to $B$ are exactly possible assignments of vertices of $F$ to the parts of (sufficiently large) blowups of $B$ that give induced copies of~$F$. For graphs $F$ and $G$, a function $f:V(F)\to V(G)$ is an \emph{embedding} of $F$ into $G$ (written as $f:F\hookrightarrow G$) if $f$ is injective and preserves both edges and non-edges; that is, $f$ gives an isomorphism of $F$ on its image. 

To avoid confusion, let us repeat that our definition of homomorphism requires that both edges and non-edges are preserved. (In the rare cases when we have to consider maps that are required to preserve edges only, we will use the term \emph{non-induced homomorphism}.) Thus, for example, an embedding can be defined as an injective homomorphism.

The \emph{edit distance} $\dedit(G,H)$ between two graphs $G$ and $H$ of the same order is the minimum value of $|E(G)\bigtriangleup f(E(H))|$ over all bijections $f:V(H)\to V(G)$; in other words, it is the smallest number of \emph{edits} (changes in adjacency) we have to do in one graph to make it isomorphic to the other. The distance from a graph $G$ to a graph family $\C G$ is 
\[
\dedit(G,\C G):=\min\{\dedit(G,H)\mid H\in\C G,\ v(H)=v(G)\}.
\] 
We will be mostly interested in the case when $\C G=\blow{B}{}$ is the family of all blowups of $B$; thus $\dedit(G,\blow{B}{})$ is the smallest number of edits in $G$ needed to make it a blowup of~$B$.

For $m\in \I N$, let $\C F_m^0$ be the family of graphs (without loops) of order $m$ consisting of one representative from each isomorphism class. For graphs $F$ and $G$ with $\kappa$ and $n$ vertices respectively, let $\P(F,G)$ be the number of $\kappa$-subsets $X\subseteq V(G)$ that induce a subgraph isomorphic to $F$ in~$G$. 

Suppose that $\kappa\le n$. Then we let $\p(F,G):=\P(F,G)/\binom{n}{\kappa}$, calling it the \emph{(induced) density} of $F$ in~$G$.
Occasionally, it will be more convenient to work with the \emph{embedding density} $\tinj(F,G)$ which is defined as the probability that a random injective map $V(F)\to V(G)$ is an embedding, that is, preserves both edges and non-edges. Thus, informally speaking, we look at vertex labelled copies of $F$ in $G$. One can easily transfer between these two densities using that
 \begin{equation}\label{eq:HomDensity}
 \tinj(F,G)= \frac{|\aut(F)|}{\kappa!}\, \p(F,G),
 \end{equation}
 where $\aut(F)$ denotes the automorphism group of $F$.
 
We call a sequence of growing bipartite graphs $(G_n)_{n\in\I N}$ with almost equal parts  \emph{$c$-quasirandom} if for every bipartite  graph $F$ 
%with parts $U_0\cup U_1=V(F)$ 
the \emph{bipartite non-induced homomorphism density} $t_{\mathrm{bip}}(F,G)$ of $F$ in $G_n$ (which is the probability that a random part-preserving map $V(F)\to V(G_n)$ sends edges of $F$ to edges of $G$) is $c^{|E(F)|}+o(1)$ as $n\to\infty$. Note that this is the value we observe in a typical $c$-random subgraph of $\MultiPartite{n,n}$. As it can be shown by an easy adaptation of the classical proof of Chung, Graham and Wilson~\cite{ChungGrahamWilson89} (with details spelled in e.g.~\cite[Lemma 14]{CooleyKangPikhurko22}), it is enough to check this property only when $F$ is  the edge and the 4-cycle.

Next, let us introduce some notation that will allow to treat the problems studied in this paper in a uniform way. 
Given an integer $\kappa\ge 2$ and a function
$\gamma:\C F_\kappa^0\to\I R$, we consider the following function on graphs:
\[
 \Obj[\gamma]{G}:=\sum_{F\in\C F_\kappa^0} \gamma(F)\P(F,G),\quad \mbox{for a graph $G$},
 \]
 and its density version
\[
\obj[\gamma]{G}:=\frac{\Obj[\gamma]{G}}{\binom{v(G)}{\kappa}}=\sum_{F\in\C F_\kappa^0} \gamma(F)\p(F,G),\quad \mbox{for a graph $G$ with $v(G)\ge \kappa$}.
\]
 We consider the corresponding extremal problem where we maximise these objective functions on order-$n$ graphs, namely we are interested in
\[ 
 \Obj[\gamma]{n}:=\max\{\Obj[\gamma]{G}\mid v(G)=n\},\quad \mbox{for $n\in\I N$},
 \] 
 and its density version $\obj[\gamma]{n}:=
 %\max\{\obj{G}\mid v(G)=n\}=
 \Obj[\gamma]{n}/\binom{n}{\kappa}$ for $n\ge \kappa$. 
 Of course, these two functions are equivalent and interchangeable. 
 %While formulas with $\obj[\gamma]{}$ are usually cleaner, 
 We tend to use $\obj[\gamma]{\cdot}$ but switch to $\Obj[\gamma]{\cdot}$ when the involved quantity seems to be better to write or understand as some counting. %and the normalised version $\obj[\gamma]{G}$ if switching to .
 Define $\obj[\gamma]{}:=\lim_{n\to\infty} \obj[\gamma]{n}$, where the existence of the limit follows from e.g.~\PST{Lemma~2.2}.

Let us observe that the problems studied here can be represented this way. For the $(\kappa,\ell)$-edge-inducibility problem, we use the same $\kappa$ and let $\gamma$ be 1 on every graph with exactly $\ell$ edges and be 0 otherwise. For the $F$-inducibility problem, we let $\kappa:=v(F)$ and 
let $\gamma$ take 
the value 1 on $F$ (or, more precisely, the unique graph in $\C F_\kappa^0$ isomorphic to $F$) and $0$ on any other graph. 
\hide{For the semi-inducibility of $H$ with $\kappa$ vertices, we define the value of $\gamma$ on $F\in\C F_\kappa^m$ to be $\Obj[H]{F}$. Then $\Obj[H]{G}=\Obj[\gamma]{G}$ for every graph $G$. However, we use $1/n^{(\kappa)}$ as the normalisation factor for $\obj[H]{n}$ since it has natural interpretation as a probabily of a random injection being an embedding of $H$ into $G$.  Since the default normalisation factor for $\obj[\gamma]{G}$ is $\binom {v(G)}{\kappa}^{-1}$, it holds for every graph $G$ with at least $\kappa$ vertices that  $\obj[\gamma]{G}=\kappa!\, \obj[H]{G}$.}

%Let a function $\gamma:\C F_\kappa^0\to\I R$ an $m$-vertex pattern $B$ be given. 
 For a graph $G$ and a vertex $u$ of $G$, let 
$\Obj[\gamma]{G,u}$ be the sum over all $\kappa$-subsets $X\subseteq V(G)$ that contain $u$ of the value of $\gamma$ evaluated at the induced subgraph $G[X]$ (more precisely, we take $\gamma(F)$ for the unique $F\in\C F_\kappa^0$ isomorphic to $G[X]$). For example, it holds for any graph $G$ that $\Obj[\gamma]{G}=(1/{\kappa})\sum_{u\in V(G)}\Obj[\gamma]{G,u}$ and thus $(1/\kappa)\Obj[\gamma]{G,u}$ can be considered as the contribution of a vertex $u$ to the global value $\Obj[\gamma]{G}$.
For $v(G)\ge\kappa$, we also define $\obj[\gamma]{G,u}:=\Obj[\gamma]{G,u}/{n-1\choose \kappa-1}$ to be its normalised version. For this normalisation, it is holds that $\obj[\gamma]{G}$ is the average of $\obj[\gamma]{G,u}$ over $u\in V(G)$:
 \begin{equation}\label{eq:AverageObjU}
 \obj[\gamma]{G}=\frac1{v(G)} \sum_{u\in V(G)} \obj[\gamma]{G,u},\quad \mbox{any graph $G$ with $v(G)\ge \kappa$}.
 \end{equation}

The \emph{$(m-1)$-dimensional simplex} is
\begin{equation}\label{eq:Simplex}
 \I S_m:=\left\{(x_0,\dots,x_{m-1})\in \I R^m\mid x_0+\dots+x_{m-1}=1\mbox{ and } \forall i\in [m]\ x_i\ge 0\right\}.
 \end{equation}
For $\V x=(x_0,\dots,x_{m-1})\in \I S_m$ and a pattern $B$ on $[m]$, let
$\obj[\gamma]{\blow{B}{\V x}}$ be the limit as $n\to\infty$ of $\obj[\gamma]{\blow{B}{V_0,\ldots,V_{m-1}}}$, where $|V_i|=(x_i+o(1))n$ for $i\in [m]$. This is a continuous function on $\I S_m$ (in fact, a polynomial). Let 
\[
\obj[\gamma]{\blow{B}{}}=\sup\{\obj[\gamma]{\blow{B}{\V x}}\mid \V x\in\I S_m\}.
\]
By the compactness of $\I S_m$ the supremum is attained by at least one $\V x\in\I S_m$; such vectors will be called \emph{$(\gamma,B)$-maximisers}%(or just \emph{optimal})
.
We call the pattern $B$ \emph{$\gamma$-minimal} if, for every pattern $B'$ (of order $m-1$) obtained from $B$ by removing one vertex, it holds that  $\obj[\gamma]{\blow{B'}{}}<\obj[\gamma]{\blow{B}{}}$. By compactness and continuity,  this holds if and only if the (closed) set of $(\gamma,B)$-optimal vectors is disjoint from the boundary of the simplex~$\I S_m$.
The pattern $B$ is called \emph{$\gamma$-optimal} if $\obj[\gamma]{}=\obj[\gamma]{\blow{B}{}}$, that is, we can attain the asymptotically optimal constant $\obj[\gamma]{}$ by some blowups of~$B$.  
If $\gamma$ and/or $B$ is understood, we may omit them from the above notation.

Asymptotic notation, such as $o(1)$, is taken with respect to $n\to\infty$ (where $n$ is usually the order of the unknown graph $G$); the constants hidden in it may depends on $\kappa$ and $\gamma$ but not on any other parameters. We call a sequence of graphs $(G_n)_{n\in\I N}$ with strictly increasing orders \emph{almost $\gamma$-extremal} (resp.\ \emph{almost $c$-regular}) if $\obj[\gamma]{G_n}=\obj[\gamma]{}+o(1)$ as $n\to\infty$ (resp.\ for every $\e>0$ there is $n_0$ such that, for every $n\ge n_0$, at least $(1-\e)v(G)$ vertices $u$ of $G$ satisfy $\deg_G(u)=(c\pm\e)v(G_n)$).

\subsection{Flag algebras}\label{se:FA}

Since the flag algebra approach of Razborov~\cite{Razborov07}
is well established by now (and is described in detail in e.g.~\cite{Razborov10,BaberTalbot11,SFS16,GilboaGlebovHefetzLinialMorgenstein22,JeongParkYang}),
here we just give a bare minimum of definitions needed to define what a flag algebra certificate contains.

Let $\tau$ be a \emph{type}, that is, is a graph with vertex set $[q]$ for some $q\in \I N$; we view $\tau$ as having all its $q$ vertices labelled.
A \emph{$\tau$-flag} is a pair $(F,f)$ where $F$ is a graph and $f:[q]\to V(F)$ is an embedding of $\tau$ into $F$ (that is, an injection that preserves edges and non-edges). We view a $\tau$-flag as a partially labelled graph where the labelled vertices, called \emph{roots}, induce a copy of~$\tau$. 
%When $q$ is small, we may write $(F,f)$ as $(F,f(0),\dots,f(q-1))$ instead, that is, listing the graph and its roots. 
For $\tau$-flags $(F,f)$ and $(H,h)$, let $\P((F,f),(H,h))$ be the number of $(v(F)-q)$-subsets $X$ of $V(H)\setminus h([q])$ such that the $\tau$-flag $(H[X\cup h([q])],h)$ is isomorphic to $(F,f)$, where isomorphisms between $\tau$-flags have to preserve each labelled root (in addition to preserving edges and non-edges). If $v(F)\le v(H)$, then the corresponding \emph{(flag) density} is 
\[
\p((F,f),(H,h)):=\frac{\P((F,f),(H,h))}{\binom{v(H)-q}{v(F)-q}}.
\] 
%of $(F,f)$ in $(H,h)$ is the probability for a uniformly random $(v(F)-q)$-subset $X$ of $V(G)\setminus h([q])$ that $(H[X\cup h([q])],h)$ is isomorphic to $(F,f)$, where isomorphisms between $\tau$-flags have to preserve each labelled root.  
For $s\ge q$, let $\C F_{s}^\tau$ be the set of all $\tau$-flags with $s$ vertices up to isomorphism. (This is consistent with our previous notation $\C F_s^0$ since we let $0$ denote the empty type.)
We fix, once and for all, an ordering of $\C F_s^\tau$ to be used when we have some vectors or matrices indexed by~$\C F_s^\tau$. For a $\tau$-flag $(H,h)$ with $v(H)\ge s$, let 
\begin{equation}\label{eq:TauVector}
 \V V_{(H,h)}^{\tau,s}:=\left(\P((F,f),(H,h))\right)_{(F,f)\in\C F_s^\tau}\mbox{\ \ and\ \ }\V v_{(H,h)}^{\tau,s}:=\left(\p((F,f),(H,h))\right)_{(F,f)\in\C F_s^\tau},
 \end{equation}
 be the (column) vectors listing the counts and densities  of all $s$-vertex $\tau$-sub-flags in $(H,h)$ (in the fixed ordering of $\C F_{s}^\tau$). Clearly, the entries of $\V v_{(H,h)}^{\tau,s}$ are non-negative and sum up to $1$.
 
%Suppose we want to prove 
 
Now, we can present the information that is contained in a flag algebra certificate of an upper bound $\obj[\gamma]{}\le u$ for given $\gamma:\C F_\kappa^0\to\I R$ and $u\in \I R$. The certificate lists an integer $\N\ge\kappa$ and, for each type $\tau$ on $[q]$ with $1 \le q\le \N-2$ and $q+\N$ even, a positive semi-definite matrix $X^\tau$  whose rows and columns are indexed by $\C F_s^\tau$, where $s:=(\N+q)/2$. More exactly it will be the case that, from each equivalence class $\C C$ of types under isomorphism as unlabelled graphs, we list only one representative $\tau\in \C C$ and its matrix $X^\tau$, effectively using the all-zero matrix for every other type in~$\C C$. Also, for every $F\in \C F_\N^0$, we list a real coefficient $c_F\ge 0$ (called the \emph{slack} at $F$) such that the following identity holds for every graph $G$ of order $n\to\infty$:
\begin{equation}\label{eq:FAMain}
 \sum_{F\in\C F_\N^0} (u-\obj[\gamma]{F})\P(F,G)= \sum_{{1 \le q\le \N-2\atop q\equiv N\,\mathrm{mod}\, 2}}\sum_{f:[q]\hookrightarrow V(G)}  \left(\V V_{G,f}^{\tau,s}\right)^T\!X^\tau\V V_{G,f}^{\tau,s}+ \sum_{F\in\C F_\N^0} c_F\P(F,G)+ O(n^{\N-1}),
 \end{equation}
 where, in the inner sum,  $s:=(\N+q)/2$ and $\tau:=([q],f^{-1}(E(G)))$ is the graph on $[q]$ such that the injection $f:[q]\hookrightarrow V(G)$ is an embedding of $\tau$ into~$G$. Note that the right-hand side of~\eqref{eq:FAMain}, apart from the error term, is non-negative (since all matrices $X^\tau$ are positive semi-definite and all slacks are non-negative by our assumptions). On the other hand, the left-hand side is exactly $(u-\obj[\gamma]{G})\binom n{\N}$. So the identity in~\eqref{eq:FAMain} indeed proves that $\obj[\gamma]{}\le u$. 
 %We will refer to the first term in the right-hand side of~\eqref{eq:FAMain} as the \emph{sum-of-squares} (or \emph{SOS}) term, since it can be represented this way by diagonalising the positive semi-definite matrices $X^\tau$.
 
Let us remark that, if we take any square matrices $X^\tau$ of the appropriate dimensions, then (for given $u\in\I R$) there is a unique choice of the slack coefficients $c_F$ (possibly negative) that makes~\eqref{eq:FAMain} hold for every~$G$. In brief, if we fix any type $\tau$ with $q\le \N-2$ vertices, $q\equiv \N\pmod 2$, and two $\tau$-flags $F_1,F_2\in \C F_s^\tau$, where $s:=(\N+q)/2$, then the sum $\sum_{f} \P(F_1,(G,f))\P(F_2,(G,f))$ over all embeddings $f$ of $\tau$ into an order-$n$ graph $G$ can be expressed, up to an additive $O(n^{\N-1})$ error term, as a linear combination (whose coefficients are independent of $n$) of counts of $\N$-vertex graphs in~$G$. Indeed, this sum counts the number of pairs of copies of $F_1$ and $F_2$ in $G$ sharing the same $q$ roots. Each such pair uses at most $2s-q=\N$ vertices while the contribution of each $\N$-set $X\subseteq V(G)$ depends only on the isomorphism class of the subgraph it induces in~$G$.
% (and we absorb the cases where the two copies intersect in more than $q$ vertices into the error term).
 It follows that, for given $\N$, the best possible upper bound $u$ that can be proved 
via% 
%an identity as in
~\eqref{eq:FAMain} is the value of an explicit (although usually very large) semi-definite program.

We will need the following lemma which states, roughly speaking, that the typical vectors of $\tau$-rooted densities in every almost extremal graph have to be close to the zero eigenspace of~$X^\tau$.
%matching a flag algebra bound
% and any $q$-vertex type $\tau$ used in the flag algebra proof, it holds for all but at most $o(n^q)$ embeddings $f$ of $\tau$ into $G$ that the corresponding vector of $\tau$-sub-flag densities in  $(G,f)$ has be close to the zero eigenspace of the 
%positive semi-definite  matrix $X^\tau$. 
%The proof can be found e.g.\ in \cite[Lemma~3.1]{PikhurkoVaughan13}.

\begin{lemma}\label{lm:Eigenspaces}
 Suppose that, for some $\N$, we have a flag algebra certificate proving that $\obj[\gamma]{}\le u$ as in~\eqref{eq:FAMain}. Let $\tau$ be any type present in the certificate, say with $V(\tau)=[q]$, and let $s:=(\N+q)/2$. (Thus $s$ is an integer and $s\ge q+1$).  Then for every $\e>0$ there are $\delta>0$ and $n_0$ such that, for every graph $G$ with $n\ge n_0$ vertices and $\obj[\gamma]{G}\ge u- \delta$, there are at most $\e n^q$  embeddings $f:\tau\hookrightarrow G$  such that
 \begin{equation}\label{eq:Eigenspaces}
 \left\| X^\tau \B v_{(G,f)}^{\tau,s}\right\|_\infty\ge \e, 
 \end{equation}
 where the vector $\B v_{(G,f)}^{\tau,s}$, as defined in~\eqref{eq:TauVector}, lists the densities of $s$-vertex $\tau$-flags in $(G,f)$.
 \end{lemma}

\bpf 
%Take any $\e>0$ and a type $\tau$ on $[q]$.  
Since the matrix $X^\tau$ is positive semi-definite, we have that $\V x^TX^\tau\V x=0$ if and only if $X^\tau\V x$ is the zero vector. By the compactness of $\I S_m$ and the continuity of  the function that maps $\V x\in\I S_m$ to $\V x^T X^\tau\V x$ there is $c>0$ such that every vector $\V x$  with $\|\V x\|_1=1$ and $\|X^\tau \V x\|_\infty\ge \e$ satisfies $\V x X^\tau\V x> c$. 

Let us show that, for example, $\delta:=\e c/\N!$ satisfies the lemma if $n_0$ is sufficiently large. Take any graph $G$ of sufficiently large order $n$ with $\obj[\gamma]{G}\ge \obj[\gamma]{}-\delta$. Thus the left-hand side of~\eqref{eq:FAMain} is at most $\delta {n\choose \N}$. Every embedding $f:\tau\hookrightarrow G$ for which~\eqref{eq:Eigenspaces} holds contributes at least $c {n-q\choose s-q}^2$ to the right-hand side  of~\eqref{eq:FAMain}.  By the non-negativity of all other terms in the right-hand side of~\eqref{eq:FAMain}, this identity implies that the number of such embeddings $f$ is at most $\big(\delta {n\choose \N}+O(n^{\N-1})\big) / \big(c{n-q\choose s-q}^2\big)<\e n^q$, as desired.\epf

\hide{
The following is a direct consequence of the above lemma. It can be easily strengthened in many ways to give much weaker sufficient conditions for almost regularity; however, the stated version suffices for all applications needed in this paper.

\begin{corollary}\label{co:reg}
Suppose that, for some even integer $\N$, we have a flag algebra certificate proving that $\obj[\gamma]{}\le u$ as in~\eqref{eq:FAMain}
%and there is $q\in\{2,\dots,\N-2\}$ such that for every (up to isomorphism) graph $\tau$ on $[q]$ if 
such that, for each of the 2-vertex types $\tau=(2,\{\})$ and $\tau=(2,\{01\})$, every vector $\B v$ in the kernel of $X^\tau$ satisfies that
 $$
 \sum_{F\in\C F_s^\tau} \left(\p(\tau_0,F)-\p(\tau_1,F)\right)\V v_F=0,
 $$
  where $s:=(\N+2)/2$ and, for $i\in [2]$, $\tau_i$ denotes the 3-vertex $\tau$-flag with the (unique) unlabelled vertex connected to Root $i$ but not to Root $1-i$.

Then every almost $\gamma$-extremal sequence of graphs is almost regular.\qed\end{corollary}
}

\subsection{Perfect stability}\label{se:PerfectStab}

Let $\kappa\ge 2$ and $\gamma:\C F_\kappa^0\to\I R$ be given. We define two notions of stability (namely, Erd\H os--Simonovits stability and perfect stability) and present the sufficient condition for perfect stability from~\PST{} that can be automatically verified by computer. There is one problem solved by us, namely the $(4,3)$-edge-inducibility problem with two different optimal patterns, for which only the former type of stability could hold. So we allow the pattern in our definition of Erd\H os--Simonovits stability to depend on $G$ (so that our definition applies to the $(4,3)$-problem). Namely, we call the $\obj[\gamma]{}$-problem \emph{Erd\H os--Simonovits stable} if for every $\e>0$ there are $\delta>0$ and $n_0$ such that if $G$ is a graph with $n\ge n_0$ vertices and $\obj[\gamma]{G}\ge \obj[\gamma]{}-\delta$ then there is a $\gamma$-optimal and $\gamma$-minimal pattern $B$ with $\dedit(G,\blow{B}{})\le \e {n\choose 2}$. Recall that the last inequality means that there is a partition $V(G)=V_0\cup\dots\cup V_{m-1}$ with $m:=v(B)$ such that
 \begin{equation}\label{eq:ESStab}
 \left|{\textstyle E(G)\bigtriangleup} E(\blow{B}{V_0,\dots,V_{m-1}})\right|\le \e {n\choose 2}.
 \end{equation}

Of course, if there is a unique $\gamma$-optimal and $\gamma$-minimal pattern $B$ up to isomorphism and the $(\gamma,B)$-optimal vector in $\I S_m$ is unique (up to an automorphism of the pattern $B$) then our definition implies the more common formulation of Erd\H os--Simonovits stability that any two graphs $G$ and $G'$ of the same order $n\to\infty$ with both $\obj[\gamma]{G}$ and $\obj[\gamma]{G}$ being $\obj[\gamma]{}+o(1)$ are $o(n^2)$-close to each other in the edit distance. This property is very useful as the first step towards characterizing graphs of sufficiently large order $n$ with $\obj[\gamma]{G}=\obj[\gamma]{n}$, see for example the proof of Theorem~\ref{th:PST7.1} here. This approach was pioneered by Erd\H os~\cite{Erdos67} and
Simonovits~\cite{Simonovits68}.
% for resolving some Tur\'an-type problems for graphs.

Perfect stability is a stronger property which, roughly speaking, states there is a constant $C$ so that~\eqref{eq:ESStab} holds for any function $\e(n)\ge 0$ with  $\delta:=\e/C$. Following~\PST{},  we call the $\obj[\gamma]{}$-problem \emph{perfectly $B$-stable} for a pattern $B$ if there is $C>0$ such that for every graph $G$ of order $n\ge C$ we have
 \begin{equation}\label{eq:PerfectStabDef}
 \dedit(G,\blow{B}{})\le C\,
 \left(\obj[\gamma]{n}-\obj[\gamma]{G}\right)n^2.
 \end{equation}
In particular, it follows that, for all $n\ge C$, every order-$n$ graph $G$ with $\obj[\gamma]{G}=\obj[\gamma]{n}$ is a blowup of $B$ and then the problem reduces to just maximising an explicit polynomial of degree at most $\kappa$ over (integer) part sizes summing up to~$n$.

A version of perfect stability for the Tur\'an problem for $K_t$, that is, for maximising the number of edges in a $K_t$-free graph of given order $n$ (with $B=K_{t-1}$), was proved by F{\"u}redi~\cite{Furedi15} while
%(who considered the distance to being $(t-1)$-partite instead of complete $(t-1)$-partite as we do in this paper). 
Roberts and Scott~\cite{RobertsScott18} extended this result to forbidding any
colour critical graph. The perfect stability of some hypergraph Tur\'an problems was established by Norin and
Yepremyan~\cite{NorinYepremyan17,NorinYepremyan18}. The second author, Sliacan and Tyros~\PST{Theorem 7.1} presented a sufficient condition for a flag algebra proof to give perfect stability and successfully applied it to a number of problems (including some instances of the graph inducibility problem).

Let us present a version of \PST{Theorem 7.1} that we will need here. This result is stated in \PST{} for a more general kind of optimisation where we have a hereditary graph family $\C G$ (that is, we forbid some induced subgraphs) and we maximise our objective function over $n$-vertex graphs in $\C G$ only. Since every graph is allowed for the inducibility-type questions considered here (and thus $\C G$ consists of all graphs), we do not list $\C G$ in our notation and omit those assumptions of \PST{Theorem 7.1} that vacuously hold. On the other hand, we need a generalisation where the pattern $B$ can have loops (which indicate those parts in  blowups of $B$ where we put cliques). We have to adapt some definitions from~\PST{} accordingly. Since the proof from \PST{} rather straightforwardly extends to patterns with loops, we just sketch it. In fact, if we are allowed to replace a problem by the complementary one (where all involved graphs are replaced by their complements) then there is only one new solved case (namely the  $(5,3)$-edge-inducibility problem) where we need the pattern to have loops. 

Recall that a homomorphism of a graph $F$ to a pattern $B$ is a (not necessarily injective) map that preserves both edges and non-edges; in other words it is an assignment of the vertices of $F$ to the parts of a blowup of $B$ that gives an induced copy of~$F$. Also, recall that $B$ is  $\gamma$-minimal if, for every pattern $B'$ obtained from $B$ by removing a vertex, it holds that $\obj[\gamma]{\blow{B'}{}}< \obj[\gamma]{\blow{B}{}}$.
% In other words, if we remove any vertex $i$ from $B'$ (effectively meaning that we force the $i$-th part to be empty) then the maximum possible limiting value of $\obj[\gamma]{}$ over growing blowups  strictly decreases.

\begin{theorem}[\cite{PikhurkoSliacanTyros19}*{Theorem~7.1}]\label{th:PST7.1} Given an integer $\kappa \ge 2$ and a function $\gamma:\C F_\kappa^0\to\I R$, let $\obj[\gamma]{G}:=\sum_{F\in\C F_\kappa^0} \gamma(F) \p(F,G)$ be the corresponding objective function on graphs.
 Let $B$ be a pattern (possibly with loops) on $[m]$ and let $\V a\in\I S_m$ be a vector with all entries positive. Suppose that all of the following statements hold.
 \begin{enumerate}
 \item \label{it:PST7.11} We have a flag algebra certificate $\C C$ on $\N$ vertices proving that 
$\obj[\gamma]{}\le \obj[\gamma]{\blow{B}{\V a}}$ as in~\eqref{eq:FAMain}. 
\item \label{it:PST7.12} 
We have a graph $\tau$ (without loops) with at most $\N-2$ vertices such that 
 \begin{enumerate}[(a)]
 \item\label{it:PST7.12a} if we restrict the maximisation of $\obj[\gamma]{}$ to graphs without an induced copy of $\tau$ then the limit strictly decreases, that is,
  \begin{equation}\label{eq:PST7.1a}
  \lim_{n\to\infty} \max\{\obj[\gamma]{G}\mid v(G)=n,\ \mbox{$G$ is $\tau$-free}\}< \obj[\gamma]{};
  \end{equation}
  \item\label{it:PST7.12b} up to an automorphism of $B$ (which by definition has to preserve loops and non-loops), there is a unique homomorphism from $\tau$ to $B$;
  \item\label{it:PST7.12c} every two distinct vertices $x,y\in V(B)$ have distinct neighbourhoods in $f(V(\tau))$, that is, $\Gamma_B(x)\cap f(V(\tau)) \not=\Gamma_B(y)\cap f(V(\tau))$,   
   for some (or, equivalently by Item~\ref{it:PST7.12b}, for every) homomorphism $f$ of $\tau$ to $B$.
  \end{enumerate}
\item\label{it:PST7.13} Every $F\in \C F_{\N}^0$ with $c_F=0$ admits a homomorphism to $B$.
\end{enumerate}
 Suppose further that at least one of the following three statements holds:
  \begin{enumerate}[(i)]
    \item\label{it:PST7.1i} the certificate $\C C$ contains the graph $\tau$ as a type and the corresponding matrix $X^\tau$ in $\C C$ is of co-rank 1 (that is, its kernel has dimension 1);
    \item\label{it:PST7.1ii} the pattern $B$ has no loops and if we restrict maximisation of $\obj[\gamma]{}$ to graphs without any induced copy of $B$ then the limit strictly decreases;
    \item\label{it:PST7.1iii} the pattern $B$ is $\gamma$-minimal.
    \end{enumerate}

Then the problem of maximizing $\obj[\gamma]{G}$ over $n$-vertex graphs $G$ is perfectly $B$-stable. Moreover, if Item~\ref{it:PST7.1i} holds (in addition to Items~\ref{it:PST7.11}--\ref{it:PST7.13}) then $\V a$ is the unique vector  $\V x\in \I S_m$ that maximises $\obj[\gamma]{\blow{B}{\V x}}$.
\end{theorem}

Let us remark that the limit in the left-hand side of~\eqref{eq:PST7.1a} exists by an easy double-counting argument, see e.g.\ \PST{Lemma~2.2}. Also, note that we allow non-injective homomorphisms in Item~\ref{it:PST7.12b} and such maps may be potentially required (e.g.\ we may need to force that a specific vertex $x$ of $\tau$ is mapped to a loop, which can be done by adding a clone $x'$ of $x$ and making them adjacent).

\bpf[Sketch of proof of Theorem~\ref{th:PST7.1}] Suppose on the contrary that perfect stability does not hold. Thus for every $C$ there is a graph $G$ that violates~\eqref{eq:PerfectStabDef}. We let $C\to\infty$; then the order $n:=v(G)$ of $G$ goes to infinity and $\obj[\gamma]{G}=\obj[\gamma]{}+o(1)$.

Let us say that
%, in the above notation, 
the problem is \emph{robustly $B$-stable} if there is a constant $C_1>0$ such that for every graph $H$ of order $n\ge C_1$ it holds that
 \[
 \dedit(H,\blow{B}{})\le C_1 \max\left\{n,(\obj[\gamma]{n}-\obj[\gamma]{H})n^2\right\}.
 \]
 Informally speaking, this is the same as the definition of perfect $B$-stability except we do not stipulate anything about order-$n$ graphs $H$ with $\obj[\gamma]{H}=\obj[\gamma]{n}+O(1/n)$.

Let us show that the $\obj[\gamma]{}$-problem is robustly $B$-stable, following \PST{Theorem 4.1}.  Let $u:=\obj[\gamma]{}$. By \PST{Lemma~2.2} (which is based on the standard blowup trick) we have that $\obj[\gamma]{n}=u+O(1/n)$. It follows that our hypothetical counterexample $G$ satisfies that 
 \begin{equation}\label{eq:FSlack}
 \p(F,G)=O(\max\{1/n, \obj[\gamma]{n}-\obj[\gamma]{G}\}),\quad \mbox{for every $F\in\C F_\N^0$ with $c_F>0$}.
 \end{equation}
 Thus, by Assumption~\ref{it:PST7.13}, every $\N$-vertex graph $F$ not admitting a homomorphism to $B$ satisfies~\eqref{eq:FSlack}.

Furthermore, by Assumption~\ref{it:PST7.12a} and the supersaturation argument of Erd\H os and Simonovits~\cite{ErdosSimonovits83}, we have that $\p(\tau,G)=\Omega(1)$. Assume that the vertex set of $\tau$ is $[q]$.  Fix a homomorphism $f$ from $\tau$ to~$B$. Every copy $\tau'$ of $\tau$ in $G$, say given by an injection $g:[q]\to V(G)$, defines $m$ pairwise disjoint sets $V_0,\dots,V_{m-1}$, where $V_i$ consists of those $x\in V(G)\setminus g([q])$ such that $g^{-1}(\Gamma_G(x))=f^{-1}(\Gamma_B(i))$. Informally speaking, we put $x$ into $V_i$ if the $G$-adjacencies of $x$ to $\tau'$ are the same as the $B$-adjacencies of $i$ to $\tau$. Distribute the vertices of the remainder $V(G)\setminus(V_0\cup\dots,V_{m-1})$ arbitrarily into the parts, e.g.\ let all be assigned to~$V_0$. 
Assumptions~\ref{it:PST7.12b}, \ref{it:PST7.12c} and~\ref{it:PST7.13} of the theorem imply by $q\le \N-2$ that, for every \emph{wrong} pair $xy$, that is,   $xy\in E(G)\bigtriangleup E(\blow{B}{V_0,\dots,V_{m-1}})$, every $\N$-set $Y\supseteq \{x,y\}\cup g([q])$ spans a subgraph in $G$ not admitting a homomorphism to $B$. 
It follows from~\eqref{eq:FSlack} that the expected number  of  wrong pairs over a uniformly random induced copy of $\tau$ in $G$ is at most $O(\max\{n,(\obj[\gamma]{n}-\obj[\gamma]{G})n^2\})$. 
If we change the adjacency of every wrong pair, we obtain a blowup of $B$; thus  robust stability holds.

Let us show that the pattern $B$ from the statement of Theorem~\ref{th:PST7.1} is $\gamma$-minimal. This is assumed in  Item~\ref{it:PST7.1iii} so we have to consider the other two cases.

First, suppose that Item~\ref{it:PST7.1ii} holds. Suppose on the contrary that some vector $\B x\in\I S_{m-1}$ and a pattern $B'$ obtained by removing a vertex from $B$  satisfy $\obj[\gamma]{\blow{B'}{\B x}}=\obj[\gamma]{\blow{B}{\B a}}$, 
that is, we can asymptotically attain the maximum value $\obj[\gamma]{}$ using blowups of~$B'$. Then, by Item~\ref{it:PST7.1ii}, some sufficiently large blowups of $B'$ contain a copy of $B$.  
This means that $B$ admits a homomorphism $h$ to $B'$. But then some two vertices of $B$ that are sent by $h$ to the same vertex of $B'$ contradict Assumption~\ref{it:PST7.12c}.
%\op{As referee pointed out, the last step is wrong: eg distinct non-adjacent $x,x'\in V(B^\circ)$, with $x$ loop in $B$ are both sent to non-loop $y\in V(B')$ and image of $\tau$ in $B$ contains $x$ and thus distinguishes $x$ from $x'$ :-((( }

So consider the case that Item~\ref{it:PST7.1i} holds. 
Let us show that $\V a$ is the unique vector in $\I S_m$ maximising $\obj[\gamma]{\blow{B}{\V  a}}$, following \PST{Lemma~6.2}.
Let $H$ be any $n$-vertex blowup of $B$ with $\Obj[\gamma]{H}=\Obj[\gamma]{n}+o(n^{\kappa})$. When we evaluate the identity~\eqref{eq:FAMain} on $H$ then the left-hand side is $o(n^{\N})$. So, each term in the right-hand side is also $o(n^{\N})$ by the non-negativity of the terms. In particular, this applies to the term obtained by  summing over all embeddings $f$ of $\tau$ into $H$ of $(\V V_{H,f}^{\tau,s})^T X^\tau \V V_{H,f}^{\tau,s}$, where $s:=(\N+q)/2$ and $\V V_{H,f}^{\tau,s}$ is the vector of the counts of $s$-vertex $\tau$-flags in $(H,f)$, as defined in~\eqref{eq:TauVector}. Recall that we have $\Omega(n^{q})$ embeddings $f$ by Item~\ref{it:PST7.12a} and supersaturation. Each such embedding comes from some homomorphism $g$ of $\tau$ to~$B$. There are only a bounded number of possible functions $g$ (at most $m^q$) so fix the most frequent~$g$. All maps $f$ corresponding to this $g$ give the same (normalised) vector $\V v:=\V v_{H,f}^{\tau,s}$ and this vector $\V v$ necessarily satisfies $\V v^T X^\tau \V v=o(1)$. Thus $\V v$ is close to the 1-dimensional zero eigenspace of $X^\tau\succeq 0$. Also, 
the part sizes of $H$ are determined by $\V v$, since each part ratio can be written by Items~\ref{it:PST7.12b} and~\ref{it:PST7.12c} as the rooted density in $(H,g)$ of a certain $\tau$-flag with $q+1\le s$ vertices. Passing to the limit as $n\to\infty$, we obtain that the vector $\V a$ is the unique maximiser in $\I S_m$. Now, since the unique maximiser $\V a$ has all entries positive, the pattern is $\gamma$-minimal, as desired.
 
Thus $B$ is $\gamma$-minimal in all cases (and we also proved the uniqueness of $\V a$ when Item~\ref{it:PST7.1i} holds).

Next, we need two further definitions from \PST{}, slightly simplified as we do not have any forbidden subgraphs. We call the pattern $B$ \emph{flip-averse} if there are $\delta>0$ and $n_0$ such that for every blowup $H=\blow{B}{V_0,\dots,V_{m-1}}$ of order $n\ge n_0$ with $\obj[\gamma]{H}\ge \obj[\gamma]{}-\delta$ and any pair $xy\in \binom{V(H)}2$ it holds that  $\Obj[\gamma]{H}-\Obj[\gamma]{H\oplus xy}\ge \delta n^{\kappa -2}$, where $H\oplus xy$ is obtained from $H$ by changing the adjacency between $x$ and $y$. Informally speaking, we require that every single adjacency change in an almost optimal blowup decreases the objective function by the maximum possible amount in the order of magnitude. Of course, for every $i,j\in [m]$, it is enough to check this property only for one pair $(x,y)\in V_i\times V_j$.
Next, we call $B$ \emph{strict} if for every $\e>0$ there are $\delta>0$ and $n_0$ such that,
for every $\V x\in \I S_m$ that maximises $\obj[\gamma]{\blow{B}{\V x}}$, it holds that if $\obj[\gamma]{H',y}\ge \obj[\gamma]{}-\delta$  for a graph $H'$ obtained by attaching a new vertex $y$ to a blowup $H:=\blow{B}{V_0,\dots,V_{m-1}}$ with $n\ge n_0$ vertices, where $|V_i|=x_in\pm1$ for $i\in [m]$, then there is $i\in [m]$ such that $y$ can be made into a clone of a vertex in $V_i$ after doing at most $\e n$ adjacency edits at~$y$. Using compactness, these two definitions can be equivalently restated in the language of step graphons; for example, the latter would state that the only way to optimally attach to an optimal blowup of $B$ is to be a clone (up to a null set) of an existing vertex.

We can show that the pattern $B$ is both flip-averse and strict, following the argument from \PST{Theorem 5.13}. If the first property is violated by an $n$-vertex blowup $H=\blow{B}{V_0,\dots,V_{m-1}}$ and $x,y$ with $x\in V_i$ and $y\in V_j$ then the graph $H'$ obtained from $H$ by changing all adjacencies between $X$ and $Y$ for some disjoint $\e n$-subsets $X\subseteq V_i$ and $Y\subseteq V_j$ for small fixed $\e>0$ would violate the robust stability that we already proved.  Indeed, $H'$ can be shown to be $\Omega(\e^2n^2)$-far in the edit distance from a blowup of $B$ while $\obj[\gamma]{H}-\obj[\gamma]{H'}=o(\e^2)$ so no $\delta>0$ satisfying robust stability can exist. (Also, note that such sets $X$ and $Y$ always exists, even if $i=j$: since $B$ is $\gamma$-minimal, there is some $\e>0$ such that $x_i>3\e$ for every optimal vector $\V x\in\I S_m$.)
Similarly, if the strictness property fails then we can violate the robust stability by adding $\e n$ new copies of the same vertex~$v$.

Now, we are ready to get the final contradiction, as in \PST{Theorem 5.8}. Following~\PST{}, pick small positive constants $c_6\gg \dots\gg c_1$. Assume by taking $n$ sufficiently large that $\obj[\gamma]{G}\ge \obj[\gamma]{}-c_1/2$. By Erd\H os-Simonovits stability (which directly follows from the established robust stability), take a blowup $H:=\blow{B}{V_0,\dots,V_{m-1}}$ with $|W|\le c_2{n\choose 2}$, where $W:=E(G)\bigtriangleup E(H)$ is the set of wrong pairs. It must hold that the vector $\B b:=(|V_0|/n,\dots,|V_{m-1}|/n)$ is $c_3$-close in the $\ell_1$-distance to an optimal vector $\V a'=(a_0',\dots,a_{m-1}')$. By the minimality of $B$, we have that each $a_i'$ is at least $c_6$ (and so $b_i\ge c_6/2$). Also, for every $x\in V(G)$, we have $\obj[\gamma]{G,x}\le \obj[\gamma]{}+c_2$ (as otherwise by adding e.g.\ $c_2n$ clones of the vertex $x$ we can bring $\obj[\gamma]{G}$ well over $\obj[\gamma]{}$, a contradiction). This allows us to conclude that $|S|\le (2c_2/c_4)n$, where we let 
\[S:=\left\{x\in V(G)\mid \obj[\gamma]{G,x}\le \obj[\gamma]{}-c_4\right\}.
\]
 %calling these vertices \emph{special}. 
 Since the set $S$ is small, it follows from the strictness of $B$  that each vertex outside of $S$ has the $W$-degree less than $c_5n$ (\PST{Equation~(29)}). This in turn can be used to prove that the effect on the objective function with respect to $G$ of a single edge flip for a pair inside $V(G)\setminus S$ is within additive $c_6n^{\kappa-2}$ the same as the  effect of the flip with respect to $H$ (\PST{Claim~5.12}). Again by strictness, we can ``fix'' all vertices in $S$ to have correct adjacencies to the parts $V_i$ so that, in terms of $\Obj[\gamma]{G}$, the objective function $\Obj[\gamma]{G}$ increases by  at least $c_4 n^{\kappa-1}$ per each fixed vertex.
 
Repeating iteratively this for every vertex of $S$ and  flipping the remaining wrong pairs (which all lie inside $V(G)\setminus S$), we increase the objective function $\Obj[\gamma]{G}$ by at least $c_3(n|S|+|W|)n^{\kappa-2}$  (\PST{Equation~(33)}) while doing at most $n|S|+|W|$ changes. Thus our assumption that $G$ fails perfect stability with $C\to\infty$ leads to a contradiction.\epf

\section{Computer-generated results}\label{se:computer}

In this section we present the results whose proof was generated by computer, postponing those stability and exact structure results that do not directly follow from Theorem~\ref{th:PST7.1} 
but require some extra arguments 
to Sections~\ref{se:43}--\ref{se:LeafedC4}.

The upper bound on the (edge/graph) inducibility constant $\obj[\gamma]{}$ for the appropriate $\gamma$ in each theorem of this section is proved via a flag algebra identity~\eqref{eq:FAMain}, and we just give the value of $\N$ in the statement of the corresponding theorem. 
Also, if perfect stability is claimed then it was derived by applying Theorem~\ref{th:PST7.1}; 
%it happens that Item~\ref{it:PST7.1i} applies in all cases of this paper. (So we included Item~\ref{it:PST7.1ii} into Theorem~\ref{th:PST7.1} for future reference.)
we just indicate which of its alternatives~\ref{it:PST7.1i} or~\ref{it:PST7.1ii} was applied and the used graph~$\tau$. (In this paper, we do not use Item~\ref{it:PST7.1iii} of Theorem~\ref{th:PST7.1} but we included it for future reference.)
%It happens that the only case when Item~\ref{it:PST7.1i} does not apply is Theorem~\ref{th:semiind10} (when we need to proof the uniqueness of the maximiser $\B a$ of $\obj[\gamma]{\blow{B}{\V a}}$). Of course, when the graph $\tau$ used in Theorem~\ref{th:PST7.1} is  small, Assumption~\ref{it:PST7.12} can be easily verified without a computer; however, it is more convenient to deal with it in an automated way.
In all cases that are solved in this paper, the provided maximiser $\B a$ of $\obj[\gamma]{\blow{B}{\V a}}$ happens to be unique up to an automorphism of $B$ (and if perfect stability is proved
via Theorem~\ref{th:PST7.1} then the same value of $\N$ is used for this). Recall that the uniqueness of the maximiser $\V a$ follow automatically in the cases where Theorem~\ref{th:PST7.1}\ref{it:PST7.1i} applies; in the remaining cases we include a proof that a maximising vector in unique up to an automorphism of~$B$.
Also, our scripts verify
the stated lower bounds on $\obj[\gamma]{}$ by computing the value of $\obj[\gamma]{\blow{B}{\V a}}$;
\arxiv{nonetheless, we include these (routine) calculations in this preprint.}{these calculations can also be found in the arxiv version~\cite{BodnarPikhurko25} of this paper.}

%If perfect stability is established via Theorem~\ref{th:PST7.1}\ref{it:PST7.1ii} or \ref{it:PST7.1iii} then we include a proof that the maximising vector~$\V a$ in unique up to an automorphism of~$B$. (Recall that the uniqueness of~$\V a$ follow automatically in the cases where Theorem~\ref{th:PST7.1}\ref{it:PST7.1i} applies.)

We formed the SDPs coming from~\eqref{eq:FAMain} and then analysed the solutions returned by computer, using the SageMath package \texttt{FlagAlgebraToolbox} (commit 22b8765) of the first author. This package is still under development; a short guide on how to install it and an overview of its current functionality can be found in~\cite{Bodnar26fat}.
% the GitHub repository 
%\href{https://github.com/bodnalev/sage}{\url{https://github.com/bodnalev/sage}}. 
The scripts that we used to generate the certificates, verify them and check that Theorem~\ref{th:PST7.1} applies in each stated case 
can be found 
%in the ancillary folder of the arXiv version of this paper or 
in a separate GitHub repository \href{https://github.com/bodnalev/supplementary_files/tree/main/graph_inducibility}{\url{https://github.com/bodnalev/supplementary_files/tree/main/graph_inducibility}} while the certificates themselves can be found in its sub-folder \verb$certificates$. Alternatively, all these files can be found in the ancillary folder of the arXiv version of this paper.

The certificate for e.g.\ the $(4,3)$-edge-inducibility problem is named \verb$stats43.pickle$, while a certificate for the graph inducibility 
%(resp.\ semi-inducibility) 
starts with \verb$ind$ 
%(resp.\ \verb$semiind$) 
and so should be easy to identify. Our scripts also verify (using exact arithmetical calculations) that the matrices and slacks listed in each exact-value certicate  indeed prove the claimed upper bound via an identity as in~\eqref{eq:FAMain}.
 Alternatively, a reader can use their own verifier; a description of how the data are arranged in each \verb$pickle$ file can be found in the readme file of the same \href{https://github.com/bodnalev/supplementary_files}{ \texttt{supplementary\_files}} repository.

Let us now give details of our results on edge-inducibility (that were briefly summarised in Table~\ref{ta:known}). We find it more convenient to work with ordered $\kappa$-tuples of distinct vertices when checking the claimed lower bound (namely that $\obj{\blow{B}{\V a}}$ is at least the stated value); of course, this has no effect on the density.

\newcommand{\loopB}{L}

\newcommand{\TheoremTemplate}[4]{It holds that $\obj[#1]{}=#4$, the problem is perfectly $B$-stable with $B = #2$, and the unique maximizer for $\obj[#1]{}\left(\blow{B}{\V a}\right)$ is $\V a = \left( #3 \right)$.
}

\newcommand{\OPTheoremTemplate}[7]{It holds that $\obj[#1]{}=#4$, the problem is perfectly $B$-stable for $B := #2$ by Theorem~\ref{th:PST7.1}\ref{it:PST7.1#5} with $\N=#6$ and $\tau=#7$ and every maximizer of $\obj[#1]{}\left(\blow{B}{\V a}\right)$ is $\V a = \left( #3 \right)$ up to an automorphism of $B$.
}

\newcommand{\OPTheoremTemplateNoSymmetry}[7]{It holds that $\obj[#1]{}=#4$, the problem is perfectly $B$-stable for $B := #2$ by Theorem~\ref{th:PST7.1}\ref{it:PST7.1#5} with $\N=#6$ and $\tau=#7$ and the unique maximizer of $\obj[#1]{}\left(\blow{B}{\V a}\right)$ is $\V a = \left( #3 \right)$.
}

\begin{theorem}\label{th:stats_42}
%It holds that $\obj[4,2]{}=1/2$, the problem is perfectly $B$-stable for $B = 2K_3$ by  Theorem~\ref{th:PST7.1}\ref{it:PST7.1i} with $\N=7$ and $\tau=K_3\sqcup K_2$, and the unique maximizer of $\obj[4,2]{}\left(\blow{B}{\V a}\right)$ is $\V a = \left( 1/6,\, 1/6,\, \ldots,\, 1/6 \right)$.
\OPTheoremTemplateNoSymmetry{4, 2}{2K_3}{1/6,\, \ldots,\, 1/6}{1/2}{i}{7}{K_3\sqcup K_2}
\end{theorem}
\arxiv{\bpf We only have to check that the uniform blowup $G$ of large order $n$ of two disjoint triangles satisfies $\obj[4,2]{G}\ge 1/2+o(1)$.
% Assume that $n:=v(G)$ is divisible by $6$ and 
Let $V_0$ and $V_1$ be the two connectivity components of $G$, each having $n/2$ vertices and spanning a subgraph isomorphic to $\MultiPartite{n/6,\,n/6,\,n/6}$. A 4-vertex set $X$ spans exactly 2 edges in the following two cases.
 \begin{enumerate}
 \item
We have exactly 2 vertices of $X$ inside each part $V_i$ and they form an edge in~$G$. There are $6$ ways to split 4 vertices evenly  between parts $V_0$ and $V_1$, with this partition of $X$ happening with probability $1/2^4+o(1)$ and then we have the probability of $2/3+o(1)$ for each pair to be an edge. 
 \item  We have exactly 3 vertices of $X$ are one part $V_i$ and they span exactly 2 edges (which is equivalent to them hitting exactly two out of the three parts inside $V_i$). We have to choose $i\in [2]$ and the lonely vertex to go to $V_{1-i}$. The probability of the corresponding part assignment is $1/2^4+o(1)$. Conditioned on this, we pick the 3 remaining vertices in a copy of $\MultiPartite{n/6,\,n/6,\,n/6}$. We have 3 choices for the missed part and pick the three vertices outside it (with probability $(2/3)^3+o(1)$); however each assignment when all three vertices are in the same part is counted twice (while it should not be counted at all), so we have to subtract $2\cdot(1/3)^2+o(1)$). 
 \end{enumerate}
 Thus the total $(4,2)$-density in limit as $n\to\infty$ is
 \[
  6\cdot\frac{1}{2^4}\cdot \left(\frac23\right)^2 + 2\cdot 4\cdot \frac1{2^4}\cdot\left(3\cdot\left(\frac23\right)^3-2\cdot\frac1{3^2}\right)=\frac16+\frac13=\frac12,
 \]
 as desired.\epf}{}

\begin{theorem}\label{th:stats_43}
    %\OPTheoremTemplate{4, 2}{2K_3}{\frac16, \frac16, \ldots, \frac16}{\frac12}{3K_2}{8}{}   
    It holds that $\obj[4,3]{}=1/2$. The lower bound comes from blowups of $B=(2,\{00,11\})$ (two loops) and of $B=(2,\{10\})$ (a single edge), with $\V a=(1/2,1/2)$ being the unique maximiser of $\obj[4,3]{\blow{B}{\V a}}$ in both cases. The upper bound can be proved by flag algebras as in~\eqref{eq:FAMain} using $\N=6$.
     \end{theorem}
\bpf Given the ouput of our scripts, it remains to show that the maximiser $\V a$  is a unique maximiser. Observe that a 4-subset of a blowup of $B$ spans a $(4,3)$-subgraph if and only if it has exactly 3 vertices in one part (and 1 vertex in the other part). If $x$ denotes the fraction of vertices in the first part then the asymptotic density of such $4$-sets is \begin{equation}\label{eq:p(x)}
p(x):=4(x^3(1-x)+x(1-x)^3),
\end{equation} 
where the factor $4$ is the number of ways to choose the ``lonely" vertex.  The derivative $p'(x)=-4\,(2x-1)^3$ has the unique (triple) root at $x=1/2$, which gives the global maximum of $p$ and is the unique argument on which $p$ assumes the maximum value~$1/2$.
%The upper bound $\obj[3,4]{}\le 1/2$ is proved by flag algebras using $\N=6$.
\epf

Let us observe that perfect stability does not hold for the $(4,3)$-problem even if we adapt its definition to allow multiple possible optimal patterns~$B$. Indeed, fix small $\e>0$, let $n\to\infty$ and let $G$ be obtained from $K_{n/2}\sqcup K_{n/2}$ by adding all edges between two fixed $\e n$-subsets in different parts. It is routine to see that $G$ is $\e^2n^2$-far from a union of two cliques (and from $\MultiPartite{n/2,\,n/2}$) while $\obj[4,3]{G}=1/2+O(\e^3)$. Since $\e>0$ is arbitrary, the edit distance cannot be upper bounded by  $C(\obj[4,3]{n}-\obj[4,3]{G})n^{2}$ for an absolute constant~$C$.

Even though perfect stability does not hold for the $(4,3)$-problem, we establish Erd\H os-Simonovits stability (in Lemma~\ref{lm:43stab}) and prove using it that, for all large $n$, every extremal graph is a complete bipartite graph or a union of two cliques (in Theorem~\ref{th:43exact}). For this we will need the following lemma that can be derived by examining our flag algebra certificate proving the upper bound $\obj[4,3]{}\le 1/2$.

\begin{lemma}\label{lm:43Subgraphs}  
For every $\e>0$ there are $\delta>0$ and $n_0$ such that every graph $G$ with $n\ge n_0$ vertices and $\obj[4,3]{G}\ge 1/2-\delta$ satisfies $\obj[4,1]{G}\le \e$ and $\obj[4,5]{G}\le \e$.
 \end{lemma}   
 \bpf Our script checks that, for every graph $F$ with $\N=6$ vertices that contains a $(4,1)$- or $(4,5)$-subgraph, its slack constant $c_F$ from the certificate is positive and thus lower bounded by some absolute constant $c>0$. Thus each such $F$ has density at most $\delta/c+o(1)$ in $G$. By
 $$
  \obj[4,i]{G}=\sum_{F\in\C F_\N^0} \obj[4,i]{F} \p(F,G)+o(1),
  $$
  the lemma can be satisfied by taking e.g.\ $\delta:=0.9\,c\e/2^{\binom{\N}2}$ and then sufficiently large $n_0$.\epf

\begin{theorem}\label{th:stats_52}
 %   \TheoremTemplate{5, 2}{3K_3}{1/9,\, 1/9,\, \ldots,\, 1/9}{{280}/{3^6}}       
        \OPTheoremTemplateNoSymmetry{5, 2}{3K_3}{1/9,\, \ldots,\, 1/9}{{280}/{3^6}}{i}{8}{3K_2}    
\end{theorem}

\arxiv{\bpf Let $G$ be the disjoint union of three copies of $T_{n/9,\,n/9,\,n/9}$, with vertex sets $V_0,V_1,V_2$ respectively.  A $5$-subset $X\subseteq V(G)$ spans 2 edges in exactly the two following cases depending on the distribution of vertices per parts $V_i$.
\begin{itemize}
\item There are 3 vertices inside some set $V_i$ spanning 2 edges, while the other two vertices are outside of $V_i$ and are non-adjacent. The probability of this is $3$ (pick the index $i$), times $\binom 53$ (pick the three vertices), times $3$ (pick which vertex of these three will be in a separate part), times 6 (pick 2 ordered sub-parts of $V_i$), times $1/9^3+o(1)$ (probability that these three vertices are in the specified parts), times $(2/3)^2+o(1)$ (the probability that the other two vertices are outside $V_i$), times $2/3+o(1)$ (the non-edge density of uniform blowups of $2K_3$). 

 \item The number of vertices in the parts $V_0,V_1,V_2$ is $2,2,1$ in some order and each pair inside a part spans an edge. The probability is $5$ (choose a lonely vertex and put it anywhere, say it went to $V_0$), times $6$ (choose the pair out of the 4 remaining vertices that goes into $V_1$, while the other pair goes into $V_2$), times $6^2$ (for each pair of vertices we have 6 possible ways to put them into different sub-parts of $V_i$), times $1/9^4+o(1)$ (the probability of each of the 4 vertices going into the specified sub-parts).
\end{itemize}
 Thus the lower bound is 
 \[
 3\cdot \binom53 \cdot 3 \cdot 6\cdot \frac{1}{9^3}\cdot  \left(\frac23\right)^2\cdot \frac23 +5 \cdot 6 \cdot 6^2\cdot 
 \frac1{9^4}+o(1)=\frac{160}{729}+\frac{40}{243}+o(1)=\frac{280}{729}+o(1),
 \]
 as claimed.
 \epf}{}

\begin{theorem}\label{th:stats_53}
    It holds that $\obj[5, 3]{}=(255 \sqrt{17} - 535)/2^{10}$, the problem is perfectly $B$-stable for $B:=(3,\{01,22\})$ (that is, $B$ consists of an edge plus an isolated loop) by Theorem~\ref{th:PST7.1}\ref{it:PST7.1i} with $\N=7$ and $\tau=\MultiPartite{2,1}$ (the cherry), and the unique maximiser of $\obj[5, 3]{\blow{B}{\V a}}$ is $\V a=(\alpha_1,\alpha_1, 1-2\alpha_1)$, where $\alpha_1$ was defined in~\eqref{eq:a1}.
\end{theorem}

\arxiv{\bpf Note that a 5-set $X$ spans exactly 3 edges in a blowup $\blow{B}{V_0,V_1,V_2}$ if and only if the pair $\left(|X\cap V_2|,\{|X\cap V_0|,\,|X\cap V_1|\}\right)$ is $(3,\{2,0\})$, $(2,\{2,1\})$, or $(1,\{3,1\})$. The three corresponding contributions to the limiting density are
$$
%\obj[5,3]{}\ge 
{5\choose 3}\cdot (1-2\alpha_1)^3 \cdot 2\cdot \alpha_1^2+\binom5{2,2,1}\cdot (1-2\alpha_1)^2\cdot 2\cdot \alpha_1^3
+\binom{5}{3,1,1}\cdot (1-2\alpha_1) \cdot 2\cdot \alpha_1^4
%=\frac{255 \sqrt{17} - 535}{2^{10}},
$$
 giving $\obj[5, 3]{}\ge (255 \sqrt{17} - 535)/2^{10}$, as desired. 
 \epf}{}

\begin{theorem}\label{th:stats_54}
    It holds that $\obj[5, 4]{}=5/8$, the problem is perfectly $B$-stable with $B:=(2,\{00,11\})$ (that is, $B$ consists of two isolated loops) by Theorem~\ref{th:PST7.1}\ref{it:PST7.1i} with $\N=5$ and $\tau=(3,\{01\})$ and the unique maximiser of $\obj[5, 4]{\blow{B}{\V a}}$ is $\V a=(1/2,\,1/2)$.
 %   \TheoremTemplate{5, 4}{2\loopB}{\frac12, \frac12}{\frac58}
\end{theorem}
\arxiv{\bpf
 A 5-subset of $2K_{n/2}$ spans 4 edges if and only if it has 3 vertices in one part (and 2 in the other). The probability of this is $\binom53\cdot 2\cdot 1/2^5+o(1)=5/8+o(1)$.\epf}{}

\begin{theorem}\label{th:stats_64} 
\OPTheoremTemplateNoSymmetry{6, 4}{(3,\{00,11,22\})}{1/3,\, 1/3,\, 1/3}{{40}/{3^4}}{i}{6}{(2,\{\})}
%    \TheoremTemplate{6, 4}{3\loopB}{1/3,\, 1/3,\, 1/3}{{40}/{3^4}}
\end{theorem}
\arxiv{\bpf  A 6-set spans  4 edges in a blowup of $B$ if and only if its intersections with the three parts have sizes $3,2,1$ in some order. The limit density of such sets is $\binom 6{3,2,1}\cdot 3! \cdot 1/3^6=   40/81$.\epf}{}

\begin{theorem}\label{th:stats_65}
It holds that $\obj[6, 5]{}=(10 \sqrt{10}-28)/{9}$, the problem is perfectly $B$-stable with $B:=K_2$
%(that is, $B$ consists of 2 isolated loops)  
by Theorem~\ref{th:PST7.1}\ref{it:PST7.1ii} with $\N=7$ and $\tau=(3,\{01,02\})$, and every maximiser of $\obj[6, 5]{\blow{B}{\V a}}$ is $\V a=(\alpha_2,\, 1-\alpha_2)$ up to an automorphism of $B$, where $\alpha_2$ was defined in~\eqref{eq:a2}.
\end{theorem}
\bpf A 6-set spans $5$ edges in a blowup $\blow{B}{(x, 1-x)}$ if and only if it has 5 vertices in one part (and 1 vertex in the other). This happens with limiting probability
$$
6\cdot \left(x^5(1-x)+x(1-x)^5\right),
  $$
  whose derivative is $$
  -6\,{\left(6 \, x^{4} - 12 \, x^{3} + 14 \, x^{2} - 8 \, x + 1\right)} {\left(2 \, x - 1\right)},$$ with roots $$\frac{1}{2}  \pm \frac{1}{2} \, \sqrt{\pm \frac{2}{3} \, \sqrt{10} - \frac{5}{3}}$$ and $1/2$. The highest density is attained at and only at  $\alpha_2$ or $1-\alpha_2$, both giving the final bound $\frac{10 \sqrt{10}-28}{9}$. The two  maximising vectors are the same, up to the automorphism of $B$ that swaps the parts.
\epf

\begin{theorem}\label{th:stats_67}  \OPTheoremTemplateNoSymmetry{6, 7}{(2,\{00,11\})}{1/2,\, 1/2}{{15}/{2^5}}{i}{6}{(2,\{\})}
%   \TheoremTemplate{6, 7}{2\loopB}{1/2,\, 1/2}{{15}/{2^5}}
\end{theorem}
\arxiv{\bpf A 6-set spans 7 edges if and only if it has 4 vertices in one part (and 2 in the other). Thus we get the limiting density $\binom 64\cdot 2\cdot 1/2^6=15/32$.\epf}{}

\begin{theorem}\label{th:stats_76}
It holds that $\obj[7, 6]{}=(28\sqrt{10} - 35)/{135}$, the problem is perfectly $B$-stable with $B:=K_2$ by Theorem~\ref{th:PST7.1}\ref{it:PST7.1ii} with $\N=7$ and $\tau=(3,\{\})$,
and every maximiser of $\obj[7, 6]{\blow{B}{\V a}}$ is $\V a=(\alpha_3,\, 1-\alpha_3)$ up to an automorphism of $B$, where $\alpha_3$ was defined in~\eqref{eq:a3}.
 %   \TheoremTemplate{7, 6}{K_2}{\alpha_3, 1-\alpha_3}{\frac{28\sqrt{10} - 35}{135}}
\end{theorem}
 \bpf A $7$-set spans 6 edges if and only if it has 6 vertices in one part (and 1 in the other). Thus the limiting density in $\blow{B}{(x, 1-x)}$ is
 \[
 7\cdot \left(x^6(1-x)+x(1-x)^6\right),
 \]
 whose derivative is $$
 -7\,{\left(15 \, x^{4} - 30 \, x^{3} + 25 \, x^{2} - 10 \, x + 1\right)} {\left(2 \, x - 1\right)}$$ with roots $$\frac{1}{2} \pm \frac{1}{2} \, \sqrt{ \pm \frac{4}{15} \, \sqrt{10} - \frac{1}{3}}$$ and $1/2$. The highest density is attained at and only at $\alpha_3$ or $1-\alpha_3$, both giving the final bound $\frac{28\sqrt{10} - 35}{135}$.
The two resulting  maximising vectors are the same up to the automorphism of $B$ that swaps the two parts. \epf

\begin{theorem}\label{th:stats_79}
    \OPTheoremTemplateNoSymmetry{7, 9}{(2,\{00,11\})}{1/2,\, 1/2}{{35}/{2^6}}{i}{7}{(3,\{01\})}
\end{theorem}
\arxiv{\bpf A $7$-set spans 9 edges if and only if it has 4 vertices in one part (and 3 vertices in the other part). The limiting density for two equal parts is $\binom74\cdot 2\cdot1/2^7=35/2^6$.\epf}{}

\begin{theorem}\label{th:stats_710}
    \OPTheoremTemplate{7, 10}{K_2}{1/3,\, 2/3}{{28}/{3^4}}{ii}{7}{(3,\{\})}
\end{theorem}
\bpf
 A 7-set spans 10 edges if and only if it has 5 vertices in one part (and 2 in the other). Thus the limiting density in $\blow{B}{(x, 1-x)}$ is 
 $$
 \binom{7}{2} \left( x^5{\left(1 - x\right)}^{2}  + {x^{2}\left(1 - x\right)}^{5} \right),
 $$ 
 whose maximum is attained at and only at $1/3$ or $2/3$, the optimal roots of the derivative $$
 -21\,{\left(3 \, x - 1\right)} {\left(3 \, x - 2\right)} {\left(2 \, x - 1\right)} {\left(x - 1\right)} x,$$ both giving the highest density $28/3^4,$ as claimed. The two resulting maximising vectors are the same under the automorphism of~$B$ that swaps the parts. \epf

Next, we turn to new inducibility results. Interestingly, even though the optimal part ratios are irrational in the next two results, the final density happens to be a rational number. In fact, the limiting density of every graph is a rational number in this construction.

\begin{theorem}\label{th:ind0}
Let $F:=(5,\{01,02,03\})$,  be the 3-edge star $\MultiPartite{3,1}$ plus an isolated vertex. Let the pattern $B:=(4,\{01,23\})$ consist of two disjoint edges (i.e.\ $B:=2K_2$) and let $\beta:=(3+\sqrt{3})/12$. Then $\obj[F]{}=5/24$ and the $F$-inducibility problem is perfectly $B$-stable by Theorem~\ref{th:PST7.1}\ref{it:PST7.1ii} with $\N=7$ and $\tau=(5,\{02,03,14\})$. Moreover, every maximiser of $\obj[F]{\blow{B}{\V a}}$ is $\V a=(\beta,\beta,1/2-\beta,1/2-\beta)$ up to an automorphism of~$B$.
% $\V a=\left(({3 + \sqrt{2}})/{12},\, (3 + \sqrt{2})/{12}, (3 - \sqrt{2})/{12}, (3 - \sqrt{2})/{12}\right)$.
\end{theorem}
\bpf Let $G$ be an $\left(x \mu_0, \,  x (1-\mu_0), \,  (1-x)\mu_1, \, (1-x) (1-\mu_1)\right)$-blowup of $B$ with connectivity components $V_0$ and $V_1$. Thus we have a pair of complete bipartite graphs on $V_0$ and on $V_1$ with part ratios $\mu_0, \mu_1$ respectively, and the size ratio between $V_0$ and $V_1$ is $x$. A $5$-set spans a copy of $F$ if and only if some 4 vertices are inside one part $V_i$ and span the $3$-star $\MultiPartite{3, 1}$, while the remaining fifth vertex is in the other part~$V_{1-i}$. Note that $\p(\MultiPartite{3,1},K_{\mu n, (1-\mu) n})=h(\mu)+o(1)$, where we denote $h(\mu) := 4 \left( \mu^3 (1-\mu) + (1-\mu)^3 \mu \right)$. This gives that the total probability is $$p(F, G) = 5\left(
    x^4(1-x) h(\mu_0) + (1-x)^4 x h(\mu_1)
\right) + o(1).$$
When $x \in \{0, 1\}$, this probability is $0$. Otherwise, the probability depends strictly monotonically on $h(\mu_0)$ and $h(\mu_1)$, hence we can optimize them separately. The polynomial $h$ takes its unique maximum at $1/2$, thus assume from now on that $\mu_0 = \mu_1 = 1/2$. Thus  
$$p(F, G) = \frac{5}{2} \left( x^4(1-x) + (1-x)^4x \right)+o(1).$$ 
The maximum of this polynomial is attained at $2\beta$ and $1-2\beta$, which are two roots of the derivative $$
-\frac52\,{\left(6 \, x^{2} - 6 \, x + 1\right)} {\left(2 \, x-1\right)},
$$ giving the maximal density $\frac{5}{24}$, as claimed. Also, the cases $x=2\beta$ and $x=1-2\beta$ are the same up to any automorphism of~$B$ that swaps the two edges of $B$, implying the uniqueness of the maximiser.\epf

\begin{theorem}\label{th:4Cycle+v}
Let $F:=(4,\{01,12,23,30\})$ be the 4-cycle plus an isolated vertex. Let the pattern $B:=(4,\{01,23\})$ consist of two disjoint edges and let $\beta:=(3+\sqrt{3})/12$. Then $\obj[F]{}=5/32$ and the $F$-inducibility problem is perfectly $B$-stable by Theorem~\ref{th:PST7.1}\ref{it:PST7.1ii} with $\N=7$ and $\tau=(5,\{02,03,14\})$. Moreover, every maximiser of $\obj[F]{\blow{B}{\V a}}$ is $\V a=(\beta,\beta,1/2-\beta,1/2-\beta)$ up to an automorphism of~$B$.
\end{theorem}
\bpf As before, let $G$ be an $\left(x \mu_0, \,  x (1-\mu_0), \, (1-x)\mu_1, \, (1-x) (1-\mu_1)\right)$-blowup of $B$ with connectivity components $V_0$ and $V_1$.
%Thus we have a pair of complete bipartite graphs on $V_0, V_1$ with part ratios $\mu_0, \mu_1$ respectively, and the size ratio between $V_0$ and $V_1$ is $x$. 
A $5$-set spans a copy of $F$ if and only if some 4 vertices are inside one part $V_i$ and span the $4$-cycle $\MultiPartite{2, 2}$, while the remaining fifth vertex is in the other part~$V_{1-i}$. Note that $\p(\MultiPartite{2, 2},K_{\mu n, (1-\mu) n})= h(\mu) + o(1)$, where we denote $h(\mu) := 6\, \mu^2 (1-\mu)^2$. Thus 
$$p(F, G) = 5\left(
    x^4(1-x) h(\mu_0) + (1-x)^4 x h(\mu_1)
\right) + o(1).$$
Similarly to the argument above, when $x \in \{0, 1\}$, this probability is $0$. Otherwise, the probability depends strictly monotonically on $h(\mu_0)$ and $h(\mu_1)$ whose unique maximum is at $1/2$. Thus we get the same density as above, except with a different coefficient in front: 
$$
p(F, G) = \frac{15}{8} \left( x^4(1-x) + (1-x)^4x \right)+o(1).
$$ 
As we know, the maximum is attained at and only at $2\beta$ and $1-2\beta$,
%, which is one of the optimal roots of the derivative $${\left(6 \, x^{2} - 6 \, x + 1\right)} {\left(1 - 2 \, x\right)},$$ 
giving the maximal density $\frac{5}{32}$ and the uniqueness of the maximiser up to an automorphism of~$B$.\epf

\hide{
Let $G$ be an $\V a$-blowup of $B$ of order $n$, with connectivity components $V_0$ and $V_1$. A $5$-set spans a copy of $F$ if and only if some 4 vertices are inside one part $V_i$ and span the $4$-cycle while the remaining fifth vertex is in the other part~$V_{1-i}$. Now the calculation is the same as in the proof of Theorem~\ref{th:ind0} except we use  $\p(\MultiPartite{2,2},K_{n,n})=3/8+o(1)$ obtaining that $\p(F,G)=5/32+o(1)$.}
\hide{
Note that $\p(\MultiPartite{2,2},K_{n,n})=3/8+o(1)$. Thus
 \[
 \p(F,G)= 5 \cdot \left((2\beta)^4(1-2\beta) + 2\beta\cdot (1-2\beta)^4\right)\cdot \frac38+o(1)=
 %5\cdot 4 \cdot 2\cdot 2 \left(\beta^4\cdot (1/2-\beta)+(1/2-\beta)^2 \cdot \beta\right)=
 \frac{5}{32},
 \]
The calculation for $t(F,G)$ is the same as for 3-star plus an isolated vertex, both have the same hom density. However, there are $8$ automorphisms of $C_4$: map any vertex to any other vertex (4 choices) plus two choices of the order on the cycle. Thus
$$
 i(F,G)+o(1)=\frac{5!}{8}\cdot 60(a^4b+ab^4),
 $$
 with the maximum of $5/32$.}

Our last graph inducibility result deals with  the 4-cycle plus a pendant edge. Here the inducibility constant is attained by a sequence of bipartite $5/6$-quasirandom graphs. In Theorem~\ref{th:4Cycle+eStab}, we will prove that all almost extremal order-$n$ graphs are, within $o(n^2)$ edits, of this form. For this structural result, we need some additional information from the flag algebra proof of the upper bound, which we collect here.

\begin{theorem}\label{th:4Cycle+e}
Let $F:=(5,\{01,12,23,30,04\})$ be the $4$-cycle plus a pendant edge. Then 
\[
\obj[F]{}=\frac{5^6}{2^8\cdot 3^5}=\frac{15625}{62208}.
\]
 Moreover, there is a flag algebra proof  of the upper bound with $\N=5$ that satisfies the following properties.
  \begin{enumerate}[A)]
  \item\label{it:A} For the single-vertex type $1$, the zero eigenspace of the matrix $X^1$ has dimension $1$.
  \item\label{it:B} For the cherry $\nu:=(3,\{01,02\})$, the zero eigenspace of the matrix $X^{\nu}$ has dimension $1$.
  \item\label{it:C} The zero eigenspace of the matrix $X^\sigma$ for the type $\sigma:=(3,\{\})$ (that has 3 vertices and no edges) does not contain any non-zero linear combination of
  4-vertex $\sigma$-flags where the unlabelled vertex sends at most one edge to the roots.
  \item\label{it:D} The zero eigenspace of the matrix $X^\mu$ for the type $\mu:=(3,\{01\})$ (that has 3 vertices and one edge) does not contain any non-zero linear
  combination of 4-vertex $\mu$-flags where the free vertex sends at most one edge to the roots.
  \end{enumerate}
  \end{theorem}

\bpf 
\arxiv{Let us check that any sequence of bipartite balanced $5/6$-quasirandom graphs $G$ give the stated lower bound on $\obj[F]{}$. It is easier to  compute the embedding density $\tinj(F,G)$. Since $F$ is connected, there are exactly two possible assignments of its vertices to the parts of the bipartite graph~$G$ and the probability of getting a good part assignment is $2\cdot 1/2^5+o(1)$. Conditioned on a (good) part assignment, we need to determine the bipartite embedding density of $\MultiPartite{3,2}$ minus an edge. By the $5/6$-quasirandomness, this density is $(5/6)^5(1/6)+o(1)$. Since $F$ has exactly two automorphisms, we have by~\eqref{eq:HomDensity} that 
$$
\p(F,G)=\frac{5!}{2}\, \tinj(F,G)= 60\cdot \frac{2}{2^5}\cdot \frac{5^5}{6^6}+o(1)=\frac{5^6}{2^8\cdot 3^5}+o(1),
$$
 as stated.
}{}%
The claimed upper bound on $\obj[F]{}$ is proved by a standard flag algebra application with $\N=5$. The provided script also verifies the additional Properties~\ref{it:A}--\ref{it:D} of the obtained certificate. \epf

\section{Stability and exact result for the $(4,3)$-problem}\label{se:43}

First, we show that Erd\H os-Simonovits stability (which was defined here to allow multiple constructions)  holds for the $(4,3)$-edge-inducibility problem. Recall that $\obj[4,3]{}= 1/2$, as proved in Theorem~\ref{th:stats_43}.
% as shown by e.g.\ the union of two cliques, each with $(1/2+o(1))n$ vertices.

\begin{lemma}\label{lm:43stab}
For every $\e>0$ there are $\delta>0$ and $n_0$ such that every graph $G$ with $n\ge n_0$ vertices and $\obj[4,3]{G}\ge \frac12-\delta$ is within edit distance at most $\e \binom n2$ from $\MultiPartite{n/2,n/2}$ or $2K_{n/2}$.\end{lemma}

\bpf Since we do not compute the dependence of $\delta$ on $\e$, we present a proof which is short but not efficient. (For example, the application of the Induced Removal Lemma can be avoided by fixing a ``typical" 4-set spanning a clique or an independent set, and  defining the two parts depending on the adjacencies of a vertex to the set.)
 
Suppose for contradiction that the lemma fails for some $\e>0$. Thus for every large integer $s$ there is a graph $G_s$ that  satisfies $v(G_s)\ge s$, $\obj[4,3]{G_s}\ge \frac12-\frac1s$ but is $\e\binom{v(G_s)}{2}$-far from the stated blowups.

Let $s\to\infty$. We know by Lemma~\ref{lm:43Subgraphs} that $\obj[4,1]{G_s}$ and $\obj[4,5]{G_s}$ are both $o(1)$. By the Induced Removal Lemma of Alon, Fischer, Krivelevich and Szegedy~\cite{AlonFischerKrivelevichSzegedy00}, we can change $o(1)$-fraction of adjacencies in $G_s$ and destroy all copies of $(4,1)$- and $(4,5)$-subgraphs. Of course, this changes any subgraph density by $o(1)$ so the new graph, call it $H$, satisfies $\obj[4,3]{H}=\frac12+o(1)$ as $s\to\infty$. Let $n:=v(H)$.

By Ramsey's theorem, $H$ contains a clique or an independent set with at least $4$ vertices. By passing to graph complements if needed, we can assume that the former holds. Let $V_0\subseteq V(H)$ be a maximum subset spanning a clique. If some $u$ in $V_1:= V(H)\setminus V_0$ has two distinct neighbours $w_0,w_1\in V_0$ then  for any $w$ in the set $V_0\setminus\Gamma(u)$ (which is non-empty by the maximality of $V_0$) the set $\{u,w_0,w_1,w\}$ spans a $(4,5)$-subgraph in $H$. This contradiction shows that every vertex outside of $V_0$ has at most one neighbour in~$V_0$. 

Suppose next that we have two non-adjacent vertices $u,w\in V_1$. By above and since $|V_0|\ge 4$, there are distinct $w_0,w_1\in V_0\setminus( \Gamma(u)\cup \Gamma(w))$. Thus the set $\{u,w,w_0,w_1\}$ spans a $(4,1)$-subgraph in $H$, a contradiction again.

We see that, apart at most $|V_1|\le n$ crossing edges, $H$ is the union of the cliques on $V_0$ and~$V_1$. Thus, for $x:=|V_0|/n$, we have that
 \begin{equation}\label{eq:43p}
 \Obj[4,3]{H}\ge\Obj[4,3]{K_{xn}\sqcup K_{(1-x)n}}+O(n^3)=p(x)\binom n4+O(n^3),
 \end{equation}
 where $p(x):=4(x^3(1-x)+x(1-x)^3)$. As we observed in the proof of Theorem~\ref{th:stats_43}, $x=1/2$ is the unique maximiser of this polynomial (with $p(1/2)=1/2$). Thus, by $\obj[4,3]{H}=1/2+o(1)$, each part $V_i$ has $(\frac12+o(1))n$ vertices. Thus the original graph $G_s$ is $o(n^2)$-close in the edit distance to $2K_{n/2}$, a contradiction to our assumption.\epf

Our next result shows that every $\lambda_{4,3}$-extremal graph $G$ of sufficiently large order $n$ is the disjoint union of two cliques $K_{m}\sqcup K_{n-m}$ or a complete bipartite graph $\MultiPartite{m,n-m}$ for some integer $m\in [0,n]$, without a single wrong adjacency. Clearly, the number of $(4,3)$-subgraphs is exactly $f_{m}:=m{n-m\choose 3}+{m\choose 3}(n-m)$. It is possible to describe all integers $m$ that maximise it. For this, routine calculations show that %if we increase $m$ by $1$ then $\Obj[4,3]{G}$ increases by
 $$
 f_{m+1}-f_{m}=\frac{1}{6} (n-2 m-1) \left(4 m^2-4
    m n+4 m+n^2-5 n+6\right)
    $$
    This cubic in $m$ polynomial has three real roots, namely $(n-1)/2$ and $m_{\pm}:=(n-1\pm\sqrt{3n-5})/2$. Thus (for all large $n$) there are 4 optimal choices of $m$ if $3n-5$ is a square of different parity than $n$ and 2 optimal choices otherwise. For example, in the latter case, the (unique) optimal $m$ at least $n/2$ is $\lceil m_+\rceil$, the maximal integer with $f_m-f_{m-1}>0>f_{m+1}-f_m$. In any case, every optimal $m$ is $(n\pm \sqrt{3n})/2+O(1)$, which gives, after routine calculations that, for such $m$,
    \begin{equation}\label{eq:43Lower}
    \Obj[4,3]{n}\ge \Obj[4,3]{K_{m}\sqcup K_{n-m}}=\frac12\binom{n}{4}+\frac{1}{8}\, n^2+O(n),\quad\mbox{as $n\to\infty$}.
    \end{equation}
  For comparison, to see the effect of this imbalance between the part sizes, note that
  \[
  \Obj[4,3]{K_{\lfloor n/2\rfloor}\sqcup K_{\lceil n/2\rceil }}=\frac12\binom{n}{4}-\frac{1}{16}\, n^2+O(n),\quad\mbox{as $n\to\infty$}.
  \]
\hide{
\begin{verbatim}
x = m; y = n - m; 

ddiff = Factor[
  Expand[(x + 1)*Binomial[y - 1, 3] + (y - 1)*Binomial[x + 1, 3] -
    (x*Binomial[y, 3] + y*Binomial[x, 3])]]
    
% // TeXForm

Solve[ddiff == 0, m]

x =.;
x = n/2 + Sqrt[3 n]/2 + C;
f = Expand[
  x*Binomial[n - x, 3] + (n - x)*Binomial[x, 3] - (1/2)*Binomial[n, 4]]
\end{verbatim}
}

\begin{theorem}\label{th:43exact}
 There is $n_0$ such that every graph $G$ with $n\ge n_0$ vertices and $\Obj[4,3]{G}=\Obj[4,3]{n}$ is a union of two cliques or a complete bipartite graph.
 \end{theorem}

\bpf Choose positive constants in this order, each being sufficiently small depending on the previous ones:
\begin{equation}\label{eq:constants}
\e_2\gg \e_1\gg \e_0.
\end{equation}

Let $G$ be any graph with $n\to\infty$ vertices such that $\Obj[4,3]{G}=\Obj[4,3]{n}$. 

By Lemma~\ref{lm:43stab}, we can assume that $G$ is within $\e_0\binom n2$ in the edit distance from two disjoint cliques or a complete bipartite graph. Since the problem is self-complementary, we can assume that $G$ is close to the union of two cliques. Pick a partition $V_0\cup V_1=V(G)$ such that  the symmetric difference 
\[
 W:={\textstyle E(G)\bigtriangleup} \left(\binom{V_0}{2}\cup \binom{V_1}2\right)
 \]
 has the smallest possible size. We call pairs in $W$ \emph{wrong}. Of course, $|W|\le \e_0{n\choose 2}$. Let $x:=|V_0|/n$.  
 %Observe that (as $n\to\infty$) 
%\[
%\obj[4,3]{K_{xn}\sqcup K_{(1-x)n)}}=4x(1-x)^3+4x^3(1-x)+O(1/n).
%\]
 As we already established in the proof of Theorem~\ref{th:stats_43}, the polynomial $p(x)$ in the right-hand side of~\eqref{eq:43p} is at most $1/2$ 
 %for $x\in [0,1]$ 
 with equality if and only if $x=1/2$. By continuity and compactness, we can additionally assume that $|x-1/2|\le \e_1$. As before, we do not write explicit dependencies between the constants even though this can in principle be done for all steps; for example, for this step it suffices to take, say, $\e_0\le \e_1^4/100$.
 
Recall that, for $u\in V(G)$, $\Obj{G,u}$ denotes the number of $(\kappa,\ell)$-subgraphs in $G$ that contain the vertex $u$. By~\eqref{eq:43Lower} the expected value of $\Obj[4,3]{G,u}$ for a uniformly random vertex $u$ is at least $\frac12{n\choose 4}\cdot \frac4n=\frac12\binom{n-1}3$. Fix a vertex $u$ with $\Obj[4,3]{G,u}$ being at least average. Since we cannot increase the total number of $(\kappa,\ell)$-subgraphs when replacing some vertex $w$ by  a clone of $u$, it holds that
 \begin{equation}\label{eq:Min43PerVertex}
 \Obj[4,3]{G,w}\ge \frac12\binom{n-1}3 - {n-2\choose 2},\quad\mbox{for every $w\in V(G)$.}
 \end{equation}

Next, we would like to show that each vertex is incident to small number of wrong pairs. For this we first define a polynomial $Q(x,y,z)$ that gives the limiting density of newly created $(4,3)$-subgraphs when we add to $H:=K_{xn}\sqcup K_{(1-x)n}$ a new vertex $w$ with $yn$ and $zn$ neighbours in the first and the second cliques respectively. It is routine to see that a triple $X$ in $H$ makes a $(4,3)$-subgraph with $w$ if and only if $X$ spans a triangle and sends no edge to $w$, or $X$ spans exactly 1 edge and sends 2 edges to~$w$. Thus we define
 \[
 Q(x,y,z):=\left(x-y\right)^3 + \left(1-x-z\right)^3+3y^2\left(1-x-z\right)+3\left(x-y\right)z^2 + 6yz(1-y-z).
 %\left(\frac12-y\right)^3 + \left(\frac12-z\right)^3+3y^2\left(\frac12-z\right)+3\left(\frac12-y\right)z^3 + 6yz(1-y-z).
 \]

The following claim implies that every almost optimal way for a vertex to attach to $H$ is to almost follow the pattern structure (which is basically the strictness property from the proof of Theorem~\ref{th:PST7.1}, when adapted to having more than one optimal pattern~$B$).
 
\begin{claim}\label{cl:q} The maximum value of $q(y,z):=Q(1/2,y,z)$ on $[0,1/2]^2$ is $1/2$, and is attained only at $(1/2,0)$ and $(0,1/2)$.

Also, for every $x\in [1/2-\e_1,1/2+\e_1]$ and for every $(y,z)\in [0,x]\times [0,1-x]$ at $\ell_1$-distance at least $\e_2$ from $\{(0,1-x),(x,0)\}$, it holds that $Q(x,y,z)\le 1/2-\e_{1}$.
\end{claim}

\bpf
Let $(y,z)\in [0,1/2]^2$ maximise the polynomial $q$. Of course, $q(y,z)\ge q(1/2,0)=1/2$. Suppose first that $(y,z)$ is an interior point of $[0,1/2]^2$. Then it is a critical point of $q$. Routine calculations show that the difference of the partial derivatives at $(y,z)$ is
\[
 \frac{\partial q}{\partial y}(y,z)-\frac{\partial q}{\partial z}(y,z)=6(y-z)(y+z),
 \]
 and thus $y=z$. The derivative of the cubic polynomial $q(y,y)=-20 y^3+12 y^2-3
    y/{2}+{1}/{4}$ has roots $(4\pm\sqrt6)/20$ and the values of $q(y,y)$ on these points, which are $(27\pm3\sqrt6)/100$, are both strictly less than $1/2$. Thus no interior point can maximise~$q$. 

So $(y,z)$ lies on the boundary. Suppose first that $z=0$. Then the derivative of $q(y,0)= -y^3+3 y^2-{3 y}/{4}+{1}/{4}$ has two roots $1\pm\sqrt3/2$, one belonging in $[0,1/2]$ and the other being larger than $1/2$. Thus the only possible values for optimal $y$ are $0$ and $1/2$. We have to rule out the former as $q(0,0)=1/4<1/2$ (while the latter gives the maximum). Finally, by the symmetry between $y$ and $z$, it remains to consider only the case $z=1/2$. Here $q(y,1/2)=-y^3-{3 y^2}/{2}+{1}/{2}$, its derivative has roots $-1$ and $0$, and $y=0$ is the only point in $[0,1/2]$ giving the maximum value $1/2$, proving the first part of the claim.

The second part can be derived from the first one via a compactness argument. Suppose on the contrary that, for some given $\e_2>0$, it cannot be satisfied for any sufficiently small $\e_1>0$. We let $\e_1:=1/s$ with integer $s\to\infty$ and for each $s$ pick a counterexample $(x_s,y_s,z_s)$. By passing to a subsequence of $s$, we can assume that these triples converge to some $(x,y,z)$. We have that $x=1/2$, $(y,z)\in [0,1/2]^2$ and, by the continuity of
the polynomial $Q$, that $q(y,z)=1/2$. By the first part, $(y,z)$ is $(1/2,0)$ or $(0,1/2)$ contradicting our assumption that each $(y_s,z_s)$ is $\e_2$-far from both $(0,x_s)=(0,1/2+o(1))$ and $(x_s,0)=(1/2+o(1),0)$.\epf

Let us show that every vertex $w$ satisfies $\deg_W(w)\le\e_2n$, that is, is incident to at most $\e_2n$ wrong pairs. Recall that we defined $x=|V_0|/n$. Let $y:=|\Gamma(w)\cap V_0|/n$ and $z:=|\Gamma(w)\cap V_1|/n$. The difference between  $\Obj[4,3]{G,w}$ and $Q(x,y,z)\binom{n-1}{3}$ is at most $\deg_W(w)n^2 +|W|n+O(n^2)$, where the first two terms upper bound the number of 4-sets that contain $w$ and at least one wrong pair while the last term is due to approximations (such as of ${xn\choose 2}$ by $x^2n^2/2+O(n)$). Now,~\eqref{eq:Min43PerVertex} and the second part of Claim~\ref{cl:q} give that $(y,z)$ is $\e_2$-close in the $\ell_1$-distance to $(x,0)$ or $(0,1-x)$. Since moving $w$ to the other part cannot decrease $|W|$ by the choice of the partition $V_0\cup V_1$, the former alternative holds, giving the required.

Next, let us show that each part $V_i$ spans a clique, without a single missing edge. Suppose on the contrary that some distinct $u,w\in V_i$ are non-adjacent.  Let $G'$ be obtained from $G$ by making $uw$ an edge.  Of course, every 4-subset of $V(G)$ not containing the pair $uw$ induces the same subgraph in $G$ and $G'$, and so contributes the same amount to each of $\Obj[4,3]{G'}$ and $\Obj[4,3]{G}$.
So consider a $4$-set $X\subseteq V(G)$ that contain both $u$ and $w$ but no other wrong pair (which excludes at most $(\deg_W(u)+\deg_W(w))n+ |W|\le 3\e_2 n^2$ quadruples). The set $X$ cannot span a $(4,3)$-subgraph in $G$ (spanning always 1, 2 or 5 edges). However, if the remaining two vertices of $X$ lie in different parts (that is, $|X\cap V_i|=3$) then $X$ spans a $(4,3)$-subgraph in $G'$. Thus 
\[
\Obj[4,3]{G'}-\Obj[4,3]{G}\ge (|V_i|-2)\,|V_{1-i}|- 3\e_2 n^2 >0,
\]
 which contradicts the maximality of $G$.

Thus all wrong pairs go between the two parts. This gives us quite strong control: for example, the only way to get a $(4,3)$-subgraph containing a  wrong pair $uv$ is that $uv$ is the only wrong pair and there are exactly two vertices in each part~$V_i$. 

We would like to further bound possible wrong degrees. Let $u$ be a vertex of the maximum wrong degree. Let $d:=\deg_W(u)$. Assume that $d\ge 1$ as otherwise we are done. Among the $d$ wrong neighbours of $u$ (which are all in the other part), let $w$ be one of the maximum wrong degree and let $f:=\deg_W(w)$. Of course, $1\le f\le d$. Without loss of generality, assume that $u \in V_0$ and $w\in V_1$. Recall that we denote $x=|V_0|/n$. Let $y$ be $1-x=|V_1|/n$.
Let  $G'$ be the graph obtained from $G$ by removing the edge $uw$. Consider the difference $\Obj[4,3]{G}-\Obj[4,3]{G'}$ which is non-negative by the extremality of~$G$. As before, we have to analyse only those 4-sets $X$ that contain the pair $uw$. Let $u'$ and $w'$ be the other two vertices of $X$. If such $X$ spans a $(4,3)$-subgraph in $G$ then, up to re-ordering of $u'$ and $w'$, we have that $u'\in V_0\setminus\Gamma(w)$, $w'\in V_1\setminus\Gamma(u)$
and $u'w'\not\in E(G)$. Since the number of edges in $G$ connecting $V_0\setminus\Gamma(w)$ to $V_1\setminus\Gamma(u)$ is at least $|W|-2d^2$, we get
\begin{equation}\label{eq:Guw}
\Obj[4,3]{G,uw}\le (xn-f)(yn-d)-|W|+2d^2,
\end{equation} 
 where $\Obj[4,3]{G,uw}$ denotes the number of 4-sets $X\subseteq V(G)$ spanning exactly 3 edges in $G$ and containing both $u$ and $w$.
On the other hand, we have 
\begin{equation}\label{eq:G'uw}
\Obj[4,3]{G',uw}\ge \binom{xn-f}{2}+\binom{yn-d}{2}+(d-1)(xn-f)+(f-1)(yn-d)-2d^2.
%(|W|-f-d+1):
\end{equation}
 Indeed, the first two terms count those quadruples $X\ni u,w$ that have 3 vertices in some part (and 1 vertex in the other part) such that $uw$ is the only wrong edge within $X$; note that each such $X$ spans exactly $3$ edges in~$G'$.
 The next two terms count 
 %(by the definition of $d$ and $f$) 
 the number of those $4$-sets $X=\{u,w,u',w'\}$ such that $u'\in V_0$, $w'\in V_1$ and exactly one of  $uw'$ and $u'w$ is in~$W$. Each such $X$ spans exactly 3 edges in $G'$, unless $u'w'\in W$; the number of such wrong pairs $u'w'$ can be upper bounded by~$(d-1)^2+(f-1)(d-1)\le  2d^2$, using the definitions of $d$ and $f$.

Let $c:=(x-1/2)\sqrt{n}$; thus $|V_0|=n/2+c\sqrt{n}$.  Also, let $a:=|W|/n$. 
%By the extremality of $G$, we have that $0\le \Obj[4,3]{G}-\Obj[4,3]{G'}$. 
By using the bounds in~\eqref{eq:Guw} and~\eqref{eq:G'uw} (and then that  $d\ge f$) we get after routine calculations that
 \begin{eqnarray*}
 0&\le& \Obj[4,3]{G}-\Obj[4,3]{G'}\ \le\ \mbox{RHS of~\eqref{eq:Guw}} - \mbox{RHS of~\eqref{eq:G'uw}} \nonumber\\ 
 &=& 
 \frac{3-2a-4c^2-d-f}2\,n + \frac{7d^2+6df -f^2} 2 + 3(f-d)c\sqrt{n}-\frac{3d+3f}2\nonumber\\
 & \le &\frac{3-2a-4c^2-d-f}2\,n +\frac{13}2\,d^2.\label{eq:temp1}
 \end{eqnarray*}
  Thus, by $1\le f\le d\le \e_2n$, it holds that $d\le 2$ and $f= 1$. Furthermore, since $d+f\ge 2$, we have that
  \begin{equation}\label{eq:ac}
  2a+4c^2\le 1+o(1).
  \end{equation}
 
\hide{
Let $c:=(x-1/2)\sqrt{n}$; thus $|V_0|=n/2+c\sqrt{n}$. Also, let $a:=|W|/n$. By combining~\eqref{eq:Guw} and~\eqref{eq:G'uw} we get using that $d\ge f$ and $c=o(\sqrt n)$ that
 \begin{equation}\label{eq:temp1}
 0\le \Obj[4,3]{G}-\Obj[4,3]{G'}\le \frac{3-2a-4c^2-d-f}2\,n +10d^2+o(n+d\sqrt{n}).
 \end{equation}
  Thus, by $1\le f\le d\le \e_2n$, it holds that $d\le 2$ and $f= 1$. Furthermore, since $d+f\ge 2$, we have that
  \begin{equation}\label{eq:ac}
  2a+4c^2\le 1+o(1).
  \end{equation}
} 

Next,  we write an upper bound on the global function $\Obj[4,3]{G}$ via some version of inclusion exclusion and then argue that it is incompatible with the lower bound in~\eqref{eq:43Lower} and our assumption that $d\ge 1$. In order to state it, we need to give some definitions first. A set $X\subseteq V_0\cup V_1$ is an \emph{$\{i,j\}$-set} if $\{|X\cap V_0|,|X\cap V_1|\}=\{i,j\}$. Let
\begin{itemize}
\item $t_0:={xn\choose 3}yn+xn\binom {yn}3$ be the number of $\{3,1\}$-sets;
\item $t_1:=|W|(xn-1)(yn-1)$ be the number of pairs $(X,\{u,w\})$ where $X$ is a $\{2,2\}$-set, $\{u,w\}\in W$ and $u,w\in X$;
\item $t_2:=|W|\left(\binom{xn-1}2+\binom{yn-1}2\right)$ be the number of pairs $(X,\{u,w\})$ where $X$ is a $\{3,1\}$-set, $\{u,w\}\in W$ and $u,w\in X$;
\item $t_3$ be the number of triples $(X,\{u,w\},\{u',w'\})$, where $X=\{u,w,u',w'\}$ is a $\{2,2\}$-set and $\{u,w\},\{u',w'\}\in W$; 
\item $t_4$ be the number of $\{3,1\}$-sets $X$ that span two edges in $W$.
%be the number of triples $(X,\{u,w\},\{u,w'\})$ where $X$ is a $\{3,1\}$-set, $\{u,w\},\{u,w'\}\in W$, $w\not=w'$ and $u,w,w'\in X$. 
\end{itemize}
 Let us show that 
\begin{equation}\label{eq:PIE0}
\Obj[4,3]{G} \le t_0+t_1-t_2-t_3+t_4.
\end{equation}
We prove this by listing in Table~\ref{ta:PIE} the contribution of every $4$-set $X\subseteq V(G)$ to both sides of~\eqref{eq:PIE0}.
%, where $t_0,\dots,t_3$ refer to the first four expressions of the right-hand side of~\eqref{eq:PIE}.
% and we ignore terms whose total contribution is~$o(n^2)$.
Note that there are no 4-sets spanning 3 or more wrong edges (because $W$ has  maximum degree $d\le 2$ and each neighbour of a degree-2 vertex must have degree 1, as we proved earlier). The inequality in~\eqref{eq:PIE0} follows from the fact  that, in each row of Table~\ref{ta:PIE}, the sum of the entries in the last five columns is at least the entry in the $\Obj[4,3]{G}$-column. 
%; these are absorbed into the error term in~\eqref{eq:PIE}.
%(In fact, there is some slack for two incident wrong pairs and we could have strengthened~\eqref{eq:PIE} for $|W|>(1/2+o(1))n$ but the stated inequality in~\eqref{eq:PIE} suffices to get a contradiction to $d\ge 2$.)
 \begin{table}[h]
 \begin{center}
\begin{tabular}{c|c|c|c|c|c|c|c}
\hline
$\{\,|X\cap V_0|,\,|X\cap V_1|\,\}$ & $\binom{X}2\cap W$ & $\Obj[4,3]{G}$ & $t_0$ & $t_1$ & $-t_2$ & $-t_3$ & $t_4$\\
\hline
\{4,0\} & any & 0 & 0 & 0 & 0 & 0 & 0\\
$\{3,1\}$ & empty & 1 & 1 & 0 & 0 & 0  & 0\\
\{3,1\} & single edge & 0 & 1 & 0 &  $-1$ & 0 & 0\\
\{3,1\} & 2-edge path & 0 & 1 & 0 &  $-2$ & 0 & 1\\
$\{2,2\}$ & empty & 0 & 0 & 0 & 0 & 0 & 0\\
$\{2,2\}$ & single edge & 1 & 0 & 1 & 0 & 0 & 0\\
$\{2,2\}$ & 2-edge matching & 0 & 0 & 2 & 0 & $-2$ & 0\\
$\{2,2\}$ & 2-edge path & 0 & 0 & 2 & 0 & $0$ & 0\\
\hline
%\\ (1,3) & 2-edge path & 0 & $3$ & 0 & $-2$ & 0 
\end{tabular}
\caption{Proof of the inequality in~\eqref{eq:PIE0}.}\label{ta:PIE}
 \end{center}
 \end{table}

By $d\le 2$, we have that $t_3\ge |W|^2-2n$.
Also, since $uw,uw'\in W$ with $w\not=w'$ implies that each of the vertices $w,w'$ has wrong degree $1$, we have that $t_4\le |W|/2\cdot n/2+o(n^2)=(a/4+o(1))n^2$. It follows via routine calculations  that if we substitute these bounds into~\eqref{eq:PIE0} and then subtract the lower bound $\Obj[4,3]{G}\ge \binom n4/2+(1/8+o(1))n^2$ coming from~\eqref{eq:43Lower}  (and use that $|W|=an$ while $x$ and $y$ are respectively $1/2\pm c/\sqrt{n}$)  then the terms in front of $n^4$ and $n^3$ cancel each other, while the coefficient at $n^2$ is
\[
p(a,c):= 
-a^2-2 a c^2+\frac{3
   a}{4}-\frac{c^4}{3}+\frac{c^2}{2}
   -\frac{3}{16}.
\]
%  Also, recall that
%$a=|W|/n$ while  $x$ and $y$ are equal to respectively $n/2\pm c\sqrt{n}$.

Thus it holds that $0\le p(a,c)+o(1)$. 
Since $p$ is an even function of $c$, we can assume for the calculations in this paragraph that $c\ge 0$. 
The roots of the partial derivative ${\partial p}/{\partial c}$ are $0$ and $\pm c_0$, where $c_0:=\sqrt{3-12a}/2$. 
If $a\ge 1/4$ then $p$ as a function of $c\ge 0$ is decreasing and 
% Thus 
\[
 p(a,c)\le p\left(a,0\right)= -a^2+\frac{3 a}{4}-\frac{3}{16} \le \left(\frac38\right)^2+\frac{3\cdot\frac38}{4}-\frac{3}{16}= -\frac3{64}<0,
  %=-a^2+\frac{a}{2}-\frac{25}{192}\le -\frac1{16} +\frac{1}{8}-\frac{25}{192}<0,
 %=-\frac{13}{192}\le p\left(\frac14,\frac1{2\sqrt{2}}\right),
 \]
 a contradiction.
So suppose that $a\le 1/4$. Here it holds that $p(a,c)\le p(a,c_0)=2a^2-3a/4$. The maximum of this function on $[0,1/4]$ is $0$ which is attained if and only if $a=0$. Furthermore, $p(0,c)=-(4c^2-3)^2/48$ is 0 only if $c=c_0=\sqrt{3}/2$. Thus, by the standard compactness argument, $p(a,c)\ge o(1)$ implies that $a=o(1)$ and $c=\pm \sqrt{3}/2+o(1)$.
However then~\eqref{eq:ac} is violated. This contradiction shows that $d=0$ and finishes the proof of Theorem~\ref{th:43exact}.\epf

\section{Structural results for the $4$-cycle plus a pendant edge}\label{se:LeafedC4}

Here we prove an Erd\H os--Simonovits type stability result as well as some partial structural information about large extremal graphs when $F$ is the $4$-cycle plus a pendant edge.

\begin{theorem}\label{th:4Cycle+eStab} 
Let $F:=(5,\{01,12,23,30,04\})$.
Then every sequence of graphs $(G_n)_{n\in\I N}$ of growing orders with $\obj[F]{G_n}=\obj[F]{}+o(1)$ as $n\to\infty$ consists, up to $o(v(G_n)^2)$ edits in each $G_n$, of balanced bipartite $5/6$-quasirandom graphs.\end{theorem}
\bpf
Choose positive constants, each being sufficiently small depending on the previous ones:
 $$
 %\e_7\gg \e_6\gg 
 \e_5\gg \e_4\gg \e_3\gg \e_2\gg \e_1\gg \e_0.
 $$

Let $G$ be an arbitrary graph with $n\ge 1/\e_0$ vertices and $\obj[F]{G}\ge \obj[F]{}-\e_0$.
Recall that, for a type $\tau$ on $[q]$, an integer $s>q$, and an embedding  $f:\tau\hookrightarrow G$, the vector $\V v_{(G,f)}^{\tau,s}$ lists the densities of $s$-vertex $\tau$-flags in $(G,f)$ as defined in~\eqref{eq:TauVector}. We will need to compare these with the analogous densities in a bipartite random graph $R_n$ which consists of two independent $(n/2)$-sets $V_0$ and $V_1$ where each pair across is an edge with probability $5/6$ and these events are mutually independent. If two maps $f$ and $f'$ assign each element of $[q]$ to the same part $V_i$ then, by standard concentration results (e.g.\ the Azuma-Hoeffding Inequality), the corresponding density vectors are close to each other. Thus we write $R_\infty$ to refer to the limit of $R_n$ as $n\to\infty$ and, instead of $f$, we list just the set $F_0$, where  $F_i:=f^{-1}(V_i)$ for $i\in [2]$ (that is, $F_i$ consists of the vertices of $\tau$ that are mapped by $f$ into $V_i$). 

An explicit formula for $\V v_{R_\infty,A}^{\tau,s}$ can be written as follows. Namely, if at least one of $F_0:=A$ or $F_1:=[q]\setminus A$ spans at least one edge in $\tau$ then output $0$. Otherwise, starting with the graph $\tau$ on $[q]$, add $S:=[s]\setminus [q]$ to the vertex set,  take a uniform random partition $S_0\cup S_1$ of the $(s-q)$-set $S$ and make each pair in $(F_0\times S_1)\cup (F_1\times S_0)\cup (S_0\times S_1)$ an edge with probability $5/6$ with all choices being mutually independent. Now, for a $\tau$-flag $F\in\C F_{s}^\tau$, the $F$-th entry of $\V v_{R_\infty,A}^{\tau,s}$  is the probability that the obtained random $\tau$-flag on $[s]$ with roots $0,\dots,q-1$ is isomorphic to~$F$. Of course, by the symmetry between the parts, this value will not change if we replace $A$ by $[q]\setminus A$.

Since $\e_0$ is sufficiently small depending on $\e_1$, we can assume that the conclusion of Lemma~\ref{lm:Eigenspaces} holds with respect to $\e_1$ for $\tau$ being each of the four types $1,\nu,\sigma,\mu$ appearing in Properties~\ref{it:A}--\ref{it:D} of Theorem~\ref{th:4Cycle+e}. 
Since the matrix $X^\tau$ is positive semi-definite, any vector $\V x$ with entries summing to $1$ and satisfying $\|X^\tau\V x\|\le o(1)$ is $o(1)$-close, say in the supremum norm, to the zero eigenspace of $X^\tau$. 
As the bipartite quasirandom graph $R_n$ for large $n$ also satisfies Lemma~\ref{lm:Eigenspaces}, the limiting vector $\V v_{R_\infty,A}^{\tau,s}$, where $s:=(\N+q)/2$, lies in the zero eigenspace of $X^\tau\succeq 0$ for every $A\subseteq [q]$. 
If the zero eigenspace of $X^\tau$ happens to be 1-dimensional, as it is the case for $\tau=1$ and $\tau=\nu$ by~\ref{it:A} and~\ref{it:B}, then the vectors $\V v_{G,f}^{\tau,s}$ and $\V v_{R_\infty,f^{-1}(V_0)}^{\tau,s}$ must be close to each other, except for a small number of embeddings~$f$. We conclude that the rooted densities of 3-vertex $1$-flags (resp.\ 4-vertex $\nu$-flags) in $G$ are close to what we observe in~$R_n$. Since
there is only one way to assign the vertices of $1$ or $\nu$ to the parts of $R_n$ (up to swapping $V_0$ and $V_1$) we omit $f^{-1}(V_0)$ in the formulas below. 
%Also, for the type 1 on $\{0\}$, we write $u$ to mean its embedding into $G$ that maps $0$ to $u\in V(G)$.

\begin{claim}\label{cl:1}\mbox{\ } 
\begin{enumerate}[(i)]
 \item\label{it:Type1}
  For all but at most $\e_2n$ maps $f:\{0\}\to V(G)$ it holds
that $\left\|\V v_{G,f}^{1,3}-\V v_{R_\infty}^{1,3}\right\|_\infty\le \e_2$.
  \item\label{it:TypeNu} For each but at most $\e_2n^3$ embeddings $f$ of $\nu$ into $G$, it holds
that $\left\|\V v_{G,f}^{\nu,4}-\V v_{R_\infty}^{\nu,4}\right\|_\infty\le \e_2$.\qed
\end{enumerate}
\end{claim}

The first part of Claim~\ref{cl:1} implies by simple averaging over $f(0)\in V(G)$ that, for every 3-vertex (unlabelled) graph $F$, its density in $G$ is within, say, $3\e_2$ of its limit density  $\p(F,R_\infty):=\lim_{n\to\infty} \p(F,R_n)$, that is,
 \begin{equation}\label{eq:3Vertex}
 \left|\p(F,G)-\p(F,R_\infty)\right|\le 3\e_2,\quad\mbox{for every $F\in\C F_3^0$}.
 \end{equation}
 Note that $\p(\nu,R_\infty)=(3/4)\cdot (5/6)^2=25/48$, as it equals to the probability that 3 random vertices are not all in one part times the probability that the two corresponding crossing pairs are edges. Thus $G$ has many cherries by~\eqref{eq:3Vertex} so, in particular, the second part  of Claim~\ref{cl:1} is not vacuous. 
 Also, we can express the edge density in $G$ via 3-vertex densities, namely, $\p(K_2,G)=(1/3)\sum_{F\in\C F_3^0} |E(F)|\,\p(F,G)$. Thus it follows from~\eqref{eq:3Vertex} that $\p(K_2,G)$ is within $(1/3+2/3+1)\cdot 3\e_2=6\e_2$ of the value in $R_\infty$, which is $(1/2)\cdot (5/6)=5/12$.
  
If we have a cherry with edges $ab$ and $ac$ then we say that $bc$ is its \emph{base} or that the cherry \emph{is based} on $bc$.
Note that all pairs of non-adjacent vertices $u_0u_1$ in $R_n$ are of two types: if the vertices are in two different parts then there are no cherries at all based on them and otherwise, the density of cherries based on $u_0u_1$ is approximately $(1/2)\cdot (5/6)^2=25/72$. Let us show that a similar classification of non-adjacent pairs is possible in~$G$. Define $B$ (resp.\ $C$) to consist of those $u_0u_1\in E(\O G)$ such that the number of cherries based on $u_0u_1$ is $(25/72\pm \e_3)n$ (resp.\ at most $\e_3n$). 

\begin{claim}\label{cl:Pairs} $\left|E(\O G)\setminus (B\cup C)\right|\le \e_3\binom n2$.
%The number of non-adjacent pairs not  in $B\cuFor all but at most $\e_3n^2$ non-adjacent pairs $u_0u_1$ in $G$  or differs from $25n/72$ by at most $\e_3n$.
%belongs to 
%in the interval $(25/72\pm c_4)n$.
\end{claim}

\bpf[Proof of Claim.]  Call a cherry in $G$ \emph{bad} if the corresponding embedding $f:[3]\to V(G)$ fails the conclusion of Claim~\ref{cl:1}\ref{it:TypeNu}. By Claim~\ref{cl:1}, we know that there are at most $\e_2n^3$ (vertex labelled) bad cherries. Call a pair \emph{bad} if it is the base of at least $\e_3n/2$ bad cherries. Clearly, the number of bad pairs is at most $\e_2n^3/(\e_3n/2) <\e_3\binom n2$. 

Let us show that any pair $u_0u_1\in E(\O G)$ which is not bad is in $B\cup C$.  Assume that at least $\e_3n$ cherries are based on it as otherwise $u_1u_0\in C$, as desired. Of these cherries, less than $\e_3n/2$ are bad so there is a vertex $u$ that makes a good cherry with the base $u_0u_1$. If we take a uniform random vertex $w$ in $V(G)\setminus\{u_0,u_1,u\}$ then the probabilities of the possible adjacencies of $w$ to the good cherry on $\{u_0,u_1,u\}$ are each within $\e_2$ from the corresponding values for a cherry in~$R_\infty$. Therefore, the probability that $w$ is attached to both $u_0$ and $u_1$ (which can be written as the sum of two densities depending on whether $uw$ is an edge or not) is approximated within $2\e_2<\e_3$ 
%(as there are two possibilities for $wu$) 
by the analogous probability in $R_\infty$. The latter is exactly $25/72$, giving the required.\epf 
  
By~\eqref{eq:3Vertex}, $\p(\nu,G)$ is within $3\e_2$ from $\p(\nu,R_\infty)=25/48$. 
\hide{On the other hand, 
if we sum  the number of cherries based on $u_0u_1$ for all $u_0u_1\in E(\O G)$ then we get at most  
\[
|B|\cdot (25/72+\e_3)n+|C|\cdot \e_3n + \left|E(\O G)\setminus (B\cup C)\right|\cdot n, 
\]
cherries, giving by Claim~\ref{cl:Pairs} that 
\begin{equation}\label{eq:B}
|B|\ge \frac{(\frac{25}{48}-3\e_2-12\e_3){n\choose 3}}{(\frac{25}{72}+\e_3)n}\ge \left(\frac12-20\e_3\right)\binom n2.
\end{equation}
}
On the other hand, the total number of cherries $\P(\nu,G)$ can be computed as the sum of the number of cherries based on $u_0u_1$ for all $u_0u_1\in E(\O G)$, which is at least $|B|\cdot (25/72-\e_3)n$. Thus
\begin{equation}\label{eq:B} |B|\le \frac{(\frac{25}{48}+3\e_2){n\choose 3}}{(\frac{25}{72}-\e_3)n}\le \left(\frac12+2\e_3\right)\binom{n}2.
\end{equation}
%(Using a similar but simplier calculation, we can show a lower bound on $|B|$ with the same main constant $1/2$ but we will not need this.)

Let $H:=(V(G),E(G)\cup C)$. Let us show that the are at most $4\e_4 \binom n3$ triangles in $H$. There are four different types of a triangle, depending on how many of its edges are in $G$. So it is enough to bound the number of triangles in $H$ of each type by $\e_4 \binom n3$. 

The number of triangles that take all three edges from $G$ is obviously $\p(K_3,G)\binom n3< \e_4 \binom n3$. 

Next, note that every pair in $C\subseteq E(\O G)$ is in at most $\e_3 n$ cherries by definition; thus the number of triangles in $H$ with exactly two edges from $G$ is at most $\binom n2\cdot \e_3n\le \e_4\binom n3$. 

Let us turn to triangles with exactly one edge from $G$.
Call an edge $u_0u_1\in E(G)$ \emph{bad} if it is in at least $\e_3n$ triangles of~$G$. Since the total number of triangles in $G$ is at most $3\e_2\binom{n}3$ by~\eqref{eq:3Vertex}, we have at most $3\e_2\binom{n}3/(\e_3n)\le (\e_4/6) \binom n2$ bad edges and thus at most $(\e_4/2)\binom n3$ triangles in $H$ whose unique $G$-edge is bad.
Next, take any triangle $u_0u_1u_2$ in $H$ whose unique edge from $G$, say $u_0u_1\in E(G)$, is not bad. If we take a uniform random vertex $w\in V(G)\setminus\{u_0,u_1,u_2\}$ then the probability that it sends at least two edges to $u_0u_1u_2$ is at most $3\e_3n/(n-3)<4\e_3$, since each pair $u_iu_j$ (which is a non-bad edge or a pair in $C$) belongs to at most $\e_3 n$ triangles.
Thus the vector $\V v_{G,f}^{\mu,4}$ of the densities of $4$-vertex $\mu$-flags in $(G,f)$ for the embedding $f$ of $\mu$ into $G$ that sends $i$ to $u_i$ is $O(\e_3)$-close to being supported on $\mu$-flags where the free vertex sends at most one edge to the roots. By Property~\ref{it:D}, every such vector of $\ell_1$-norm 1 must be $\Omega(1)$-far from the zero eigenspace of the matrix $X^\mu$. Since $\e_3\ll \e_4$ is sufficiently small, there are at most $\e_4 \binom n3$ such maps $f$ and thus at most $(\e_4/2)\binom n3$ corresponding unordered triples $\{u_0,u_1,u_2\}$ (since each such triple gives 2 embeddings of $\mu$). We conclude that the total number of triangles in $H$ with exactly one edge from $G$ is at most $\e_4 {n\choose 3}$, as desired.

Finally, the argument  for upper bounding the number of triangles with all three pairs coming from $C$ is similar to the one from the previous paragraph, except each such triangle spans a copy of the edgeless type $\sigma$ in $G$ and we use Property~\ref{it:C}.

Thus we have shown that $H$ has at most $4\e_4\binom n3$ triangles and, by Claim~\ref{cl:Pairs} and~\eqref{eq:B},  at least $\binom n2-|B|-\e_3\binom n2\ge (1/2-21\e_3)\binom n2$ edges. The Erd\H os--Simonovits Stability Theorem~\cite{Erdos67,Simonovits68} implies that there is a partition $V_0\cup V_1$ such that $H$ (and thus $G$) has at most $\e_5\binom n2$ edges inside $V_0$ or $V_1$. Since we can choose $\e_5$ arbitrarily small, the graph $G$ (or, more precisely, any sequence of almost extremal graphs) is almost bipartite and its $(5/6)$-quasirandomness follows
(by e.g.~\cite[Lemma 14]{CooleyKangPikhurko22}) since $G$ has the correct density of edges and of all 4-vertex subgraphs containing a cherry (and thus the correct non-induced homomorphic density of 4-cycles) by the discussion after Claim~\ref{cl:1}.\epf

Some further information can be derived about the structure of large extremal graphs.

\begin{theorem}\label{th:C4+eExtr} Let $F$ be the 4-cycle with a pendant edge.
For any $\e>0$ there is $n_0$ such that any  graph $G$ of order $n\ge n_0$ with $\obj[F]{G}=\obj[F]{n}$ admits a partition $V(G)=V_0\cup V_1$ such that each $V_i$ is an independent set and each vertex of $G$ has $(5/12\pm \e)n$ neighbours in the other part.\end{theorem}

\bpf Let $G$ be any $\obj[F]{}$-extremal graph with $n\to\infty$ vertices.  For brevity and since the meaning will be clear, we will hide all negligible constants under~$o(1)$ terms.

By Theorem~\ref{th:4Cycle+eStab}, there is a partition $U_0\cup U_1=V(G)$ such that each part has $(1/2+o(1))n$ vertices and spans $o(n^2)$ edges while the bipartite graph $G[U_0,U_1]$ (or rather the corresponding sequence of graphs as $n\to\infty$) is $p$-quasirandom, were we denote $p:=5/6$. Let $V_0\cup V_1$ be a  partition of $V(G)$ that maximises the number of crossing edges. Note that this number is at least $|E(G[U_0,U_1])|=p(n/2)^2+o(n^2)$ and at most $|E(G)|=p(n/2)^2+o(n^2)$. Thus almost all edges of $G$ must go between $V_0$ and $V_1$. By the quasirandomness, we have up to swapping the parts that $|V_i\bigtriangleup U_i|=o(n)$ for each $i\in [2]$. Thus all conclusions of Theorem~\ref{th:4Cycle+eStab} also apply to the partition $V_0\cup V_1$ (with error terms that are worse but are still $o(1)$).

Here it is more convenient to count embeddings of $F$ into~$G$ (rather than $5$-subsets that span a subgraph isomorphic to $F$).
Since $F$ has only two automorphisms, we have by~\eqref{eq:HomDensity} that $t(F,G)$ is $\gamma+o(1)$, where 
\[
 \gamma:=\frac{2}{5!}\, \obj[F]{} =\frac{1}{60}\, \frac{5^6}{2^8\cdot 3^5}=\frac{5^5}{2^{10}\cdot 3^{6}}=\frac{p^5(1-p)}{2^4}.
 \]
 Thus, by~\eqref{eq:AverageObjU}, the expected number of embeddings of $F$ into $G$ that use a uniformly random vertex is $v(F)\cdot t(F,G)n^5/n=(5\gamma+o(1))n^4$.
Since every two distinct vertices of $G$ are simultaneously in at most $5\cdot 4\cdot n^3=o(n^4)$ embeddings of $F$ into $G$ and we cannot increase $t(F,G)$ by replacing a vertex by a clone of another vertex, we conclude 
 that every vertex of $G$ is in least $(5\gamma+o(1))n^4$ embeddings. 

The rest of the proof relies on rather laborious calculations, which are also included in the ancillary notebook. Namely, the notebook verifies Claim \ref{cl:5/12} below  by computing the polynomial $q$ and checking local optimums both in the interior and on the boundary. Also, the notebook calculates the rate of decrease in the induced copies of $F$ under the addition of an edge inside a part, a quantity that we have to estimate later when we argue than $V_i$ spans no edges.

Take any vertex $u$ of $G$.  If $u$ has $x_in$ neighbours in $V_i$ and $y_in$ non-neighbours in $V_i$ for the max-cut partition $V(G)=V_0\cup V_1$ then the number of embeddings $f$ of $F$ containing $u$ is by Theorem~\ref{th:4Cycle+e} (that is, by the $5/6$-quasirandomness of $G[V_0,V_1]$) equal to
\begin{align}
\sum_{i=0}^1&\Big((x_ix_{1-i}^2y_{i}p^2(1-p)^2+x_{1-i}^3y_{i}p^2(1-p))\nonumber\\
&+ 2\cdot(x_i^2y_{1-i}^2p^3(1-p))
+x_i^2y_iy_{1-i}p^3+x_iy_i y_{1-i}^2p^4
\Big)n^4
+o(n^4),\label{eq:FPerU}
%\left(\\left(\frac 12-y\right)^2\left(\frac 12-x\right) +o(1)\right)n^4
\end{align}
where the four terms under the sum count those $f$ that map $0$, one of  $1$ or $3$, $2$, and $4$ respectively to the vertex~$u$. (Recall that $E(F)=\{01,12,23,30,04\}$.) For example, let us show that the first term counts those embeddings $f$ that send $0$ to~$u$. 
Since we have $o(n^2)$ edges inside each part, there are only $o(n^4)$ embeddings that map an edge of $F$ inside some part of $G-u$; so let us exclude these. Thus, if $i\in[2]$ is the index with  $f(2)\in V_{i}$ then $f(1),f(3)\in V_{1-i}$. Since $4$ sends no edges to $\{1,2,3\}$ in the graph $F$, its image $f(4)$ can be  in both $V_i$ or  $V_{1-i}$, which will give the two summands making the first term. For example, if $f(4)\in V_i$  then we have to count the number of choices of $f(2)\in V_i\setminus \Gamma(u)$, $f(1),f(3)\in V_{1-i}\cap \Gamma(u)$ and $f(4)\in V_i\cap \Gamma(u)$ such that $f(2)$ (resp.\ $f(4)$) is adjacent to both (resp.\ none) of $f(1)$ and $f(3)$. Since the induced bipartite graph $G[V_0,V_1]$ is $p$-quasirandom, the number of choices as above is $y_in\cdot (x_{1-i}n)^2\cdot x_in\cdot p^2(1-p)^2+o(n^4)$, giving the first summand. The other cases are analogous.

Substituting $y_i=1/2-x_i+o(1)$, we get  that the expression in~\eqref{eq:FPerU} is $(q(x_0,x_1)+o(1)) n^4$, where for $a,b\in \I R$ we define $q(a,b)$ to be
\[
 \left(\frac{25ab}{54}+\frac{125}{1728}\right)(a^2+b^2)-\frac{25
   }{108}(a^3+b^3)-\frac{50 a^2
   b^2}{81}+\left(\frac{325 a
   b}{1296}+\frac{625}{10368}\right)(a+b) -\frac{625 a
   b}{1296}
   %,\quad \mbox{for $a,b\in \I R$}
   .
\]
\hide{ \begin{eqnarray*}
 q(a,b)&:=&\frac{25 a^3 b}{54}-\frac{25
   a^3}{108}-\frac{50 a^2
   b^2}{81}+\frac{325 a^2
   b}{1296}+\frac{125
   a^2}{1728}+\frac{25 a
   b^3}{54}+\frac{325 a
   b^2}{1296}-\frac{625 a
   b}{1296}\\
   &+& \frac{625
   a}{10368}-\frac{25
   b^3}{108}+\frac{125
   b^2}{1728}+\frac{625 b}{10368},\quad \mbox{for $a,b\in \I R$}.
 \end{eqnarray*}

Note that  
 
 If we let $x_0=a$, $x_1=b$, $y_0=1/2-a$ and $y_1=1/2-b$ then the main term becomes $g(a,b)$, which is
 \begin{verbatim}
 (625 a)/10368 + (125 a^2)/1728 - (25 a^3)/108 + (625 b)/10368 - (
 625 a b)/1296 + (325 a^2 b)/1296 + (25 a^3 b)/54 + (125 b^2)/1728 + (
 325 a b^2)/1296 - (50 a^2 b^2)/81 - (25 b^3)/108 + (25 a b^3)/54
 \end{verbatim}
 }

A routine calculation confirms that $q(0,5/12)=5\gamma$, which is in accordance with our counting that led to the definition of~$q$. Let us show that this is in fact  the maximum value of $q$ over $[0,1/2]^2$ and  the only maximisers of $q$ are those pairs $(a,b)$ that we expect to see in an extremal construction:

\begin{claim}\label{cl:5/12} If $(a,b)\in [0,1/2]^2$ satisfies $q(a,b)\ge 5\gamma$ then $\{a,b\}=\{0,5/12\}$.
\end{claim}

\bpf[Proof of Claim.] Suppose that $(a,b)$ is a counterexample.

Suppose first that $c:=q(a,b)-5\gamma>0$. Take $n$ sufficiently large depending on~$c$. Consider a graph $G'$ obtained from $G$ by removing a ``typical'' vertex $w$ and adding a new vertex $w'$ with $(a+o(1))n$ and $(b+o(1))n$ neighbours in the two parts. The new vertex is in at least $(5\gamma+c+o(1))n^4$ embeddings of $F$ which is by $(c+o(1))n^4>20n^3$ larger than the number of embeddings destroyed by removing $w$. Thus $\Obj[F]{G'}>\Obj[F]{G}$, a contradiction to the maximality of~$G$. We conclude that $q(a,b)=5\gamma$ and it is a maximiser of $q$ on $[0,1/2]^2$. 

Next, suppose that $(a,b)$ lies in the interior of $[0,1/2]^2$. Since $q$ as a polynomial is a differentiable function, $(a,b)$ must be a critical point. We add a new variable $z\in\I R$ and run the Buchberger Algorithm (which is a standard function in SageMath) to eliminate variables $x$ and $y$ from the following system of 3 polynomial equations: $z-q(x,y)=0$, $\frac{\partial q}{\partial x}(x,y)=0$ and 
$\frac{\partial q}{\partial y}(x,y)=0$. The algorithm produces an explicit polynomial $Q(z)$ depending only on $z$ which is in the ideal generated by these 3 polynomials, namely
\hide{
\begin{align*}
Q(z)&=z^6+\frac{69672635
   }{382205952}\,z^5-\frac{78942300769
   8125}{194775186325635072}\,
   z^4
   \\&+\frac{20
   055643492606990625}{64621732218373901647872}\, 
   z^3-\frac{28005317558330078125}{10468720619376572066955264}\, 
   z^2
   \\&-
   \frac{21065181302459716796875}{308735133530158362140991081676
   8}\,z-\frac{370573902130126953125}{9
   8795242729650675885117146136576}
\end{align*}
}
\begin{align*}
 Q(z)&=98795242729650675885117146136576\,
   z^6\\
   &+18009465447674020572152337530
   880\,
   z^5
   -40041670159933736113668096000
   0\,
   z^4\\
   &+30661545256257839254732800000\,  z^3-264291334776391680000000000\,
   z^2\\
   &-674085801678710937500000\,
   z-370573902130126953125.
\end{align*}
\hide{
This polynomial has degree 6 but its (integer) coefficients are rather large. It happens to factorise over rationals into two degree-3 polynomials; for brevity we write just these two factors:
\begin{eqnarray*}
Q_1(z)&=&4057816381784064
   z^3+807060116275200
   z^2-8690400000000
   z-11962890625,\\
   Q_2(z) &=& 24346898290704384
   z^3-404144732897280
   z^2+33844971164000
   z+30976953125.
\end{eqnarray*}
Since the triple $(x,y,z)=(a,b,5\gamma)$ is a root of each of the 3 equations, it must hold that $Q(5\gamma)=0$. However, one can verify 
that $5\gamma$ is a root of neither $Q_1$ nor $Q_2$, a contradiction. }

Since the triple $(x,y,z)=(a,b,5\gamma)$ is a root of each of the 3 equations, it must hold that $Q(5\gamma)=0$. However, one can verify using exact arithmetic (see e.g.\ the calculations in our script)
that $5\gamma$ is not a root of $Q$, a contradiction.

Thus the point $(a,b)$ lies on the boundary of $[0,1/2]^2$. Up to symmetry, there are only the following two cases to consider. First, let $b=1/2$. We have that 
\[
q(a,1/2)=\frac{25 a^2}{576}-\frac{625
   a}{10368}+\frac{25}{1296}.
\]
The derivative of this quadratic polynomial has zero at $a=25/36>1/2$, so its maximum value on $[0,1/2]$ is $q(0,1/2)=25/1296<5\gamma$, a contradiction. Finally, let $b=0$. We have
\[
q(a,0)=-\frac{25 a^3}{108}+\frac{125
   a^2}{1728}+\frac{625 a}{10368}.
\]
The derivative of this cubic polynomial has roots $-5/24$ and $5/12$ so its unique maximiser on $[0,1/2]$ is $a=5/12$ (when the value of $q$ is $5\gamma$), 
%It is easy to see that the corners $q(0, 0), q(0, 1/2), q(1/2, 1/2)$ give a strictly smaller value than $5 \gamma$, thus the unique maximiser is when $\{a,b\}=\{0,5/12\}$, 
contradicting our initial assumption.\epf

Similarly as we derived the second part of Claim~\ref{cl:q} from the first, Claim~\ref{cl:5/12} gives by a compactness argument that $q(a,b)$ is approximately $5\gamma$ only if $\{a,b\}$ is close to $\{0,5/12\}$, which in particular applies to $\{x_0,x_1\}$. Since $V_0\cup V_1$ is a max-cut partition,  we conclude that every vertex of $G$  has $o(n)$ neighbours in its part and $(5/12+o(1))n$ neighbours in the other part.

We can now show that each part  spans no edges. Suppose on the contrary that some edge $uw\in \binom{V_i}{2}$ violates this.  Let 
$$
 y:=\frac{|V_{1-i}\cap \Gamma(u)\cap \Gamma(w)|}{n/2}
 $$
 be the proportion of vertices of $V_{1-i}$ that are adjacent to both $u$ and $w$. Since the bipartite graph $G[V_0,V_1]$ is almost $5/6$-regular, we have that $2p-1\le y+o(1)\le p$.
 
 Let us compare $G$ with the graph $G'$ which is obtained from $G$ by making the pair $uw$ a non-edge. First, let us count embeddings of $F$ into $G$ that use the edge $uw$. Since the maximum degree inside each part is $o(n)$, there are only $o(n^3)$ embeddings of $F$ into $G$ that use an edge inside a part different from the pair $uw$. Every other embedding has to map $\{0,4\}$ to $\{u,w\}$, having 2 choices here. Then the images of both $1,3\in V(F)$ have to be from the set of non-neighbours of $f(4)$ intersected with the neighbours of $f(0)$,  whose size is $(p-y)n/2+o(n)$. Finally, we have to choose $f(2)$ adjacent to both $f(1)$ and $f(3)$. By quasirandomness, the total number of embeddings is at most $2\cdot ((p-y)n/2)^2\cdot n/2\cdot p^2+o(n^3)
%=((p(1-p)^2+o(1))n^3
$.

On the other hand, if we make $uw$ a non-edge then the number of new $F$-embeddings that use both $u$ and $w$ can be lower bounded as follows. Again, up to an additive error term $o(n^3)$, we can pretend that each $V_i$ is an independent set.
We look at embeddings that map the non-edge $\{1,3\}\in E(\overline F)$ to $\{u,w\}$. There are 2 choices for the bijection between $\{1,3\}$ and $\{u,w\}$. Then we send each of $0$ and $2$ to the common neighbourhood of $u$ and $w$ which has size  $y(n/2)+o(1)$ and send $4$ to any non-neighbour of $u$ and $w$ in $V_i$, having $n/2+o(n)$ choices. We have  to ensure additionally  that $f(4)$ is adjacent to $f(0)$ but not to $f(2)$, so we have to multiply the total number of choices by $p(1-p)$ by quasirandomness. Thus the number of new embeddings $f$ with $f(13)=uw$ is at least $2\cdot (y(n/2))^2 \cdot n/2\cdot p(1-p)+o(n^3)$.
Actually, this lower bound alone suffices for a contradiction (so we do not need to look at embeddings where some other non-edge of $F$ is mapped to~$uw$).

Thus if we take the difference between the number of embeddings of  $F$ into $G'$ and $G$, and normalise  this by $n^3$ then the main term is at least
$$
\frac{5y^2}{36\cdot 4}-\frac{25\, (5/6-y)^2}{36\cdot 4}
%=-\frac{5 y^2}{36}+\frac{125   y}{432}-\frac{625}{5184}
=:r(y)
$$
 This quadratic polynomial is clearly increasing with $y$ so the minimum value of $r$ on $[2/3,5/6]$ is $r(2/3)=55/5184>0$. This contradicts the maximality of $G$ and proves that each $V_i$ is an independent set, finishing the proof of Theorem~\ref{th:C4+eExtr}.
 \epf

\hide{
Theorem~\ref{th:C4+eExtr} reduces the $F$-inducibility problem for large $n$ to a bipartite version (modulo the issue of finding optimal part sizes). Solving this bipartite problem fully seems quite challenging and the exact extremal value probably depends on the number theoretic properties of~$n$. Therefore, we are content with the above partial description of extremal graphs.
}

\section{Concluding remarks}\label{se:concluding}

There are  a number of pairs $(\kappa,\ell)$ for which we could not determine the edge-inducibility constant. Table~\ref{ta:unknown} presents numerical upper bounds returned by computer (using the maximum computationally feasible value $\N =8$) and the best lower bounds that we could find, where $P^{\,\bullet\, \circ\, \circ\, \bullet}_{{n}/{4},\, {n}/{4},\, {n}/{4},\, {n}/{4}}$ denotes the uniform blowup of the pattern $(4,\{00,01,12,23,33\})$, which is the 4-vertex path with loops at its two endpoints, while $\alpha_5, \alpha_6, \alpha_7$ are the real roots of
\begin{align*}
    404\alpha_5^3 - 310\alpha_5^2 + 82\alpha_5 - 7 & = 0, \\ 
    14\alpha_6^2 - 14\alpha_6 + 3 & = 0, \\
    294\alpha_7^4 - 390\alpha_7^3 + 190\alpha_7^2 - 40\alpha_7 + 3 & = 0,
\end{align*}
satisfying $\alpha_5=\formatnum{8}{0.173017465363062}$, $\alpha_6=\formatnum{8}{0.311017763495386}$, and $\alpha_7=\formatnum{8}{0.176458265153189}$\,.
%$\alpha_5=0.173017465363062...$, $\alpha_6=0.311017763495386...$, and $\alpha_7=0.176458265153189...$.

In these cases, we were content with an example giving a lower bound just fairly close to the upper bound (without making a serious effort of finding a best possible construction). So it is possible that some of our lower bounds listed in Table~\ref{ta:unknown} can be easily improved. 

\begin{table}[h!]
\begin{center}
\begin{tabular}{c|c|c|l} %{c|c|c|c|c}
\hline
$(\kappa,\ell)$ & Construction & Lower bound & Numerical SDP value \\ %& Rounded value\\
\hline
\hline 

$(5,5)$ & 
%$P^{\bullet\, \circ\, \circ\, \bullet}_{\frac{n}{4}, \frac{n}{4}, \frac{n}{4}, \frac{n}{4}}$ 
$P^{\,\bullet\, \circ\, \circ\, \bullet}_{{n}/{4},\, {n}/{4},\, {n}/{4},\, {n}/{4}}$ 
& $\frac{45}{128}=0.3515625$ & $\formatnum{}{0.3515625030643538}$ \\

\hline

$(6,2)$ & 
%$6K_{\frac{n}{6}}$ 
$6K_{{n}/{6}}$ 
& $\frac{25}{72}=\formatnum{}{0.347222222222222222}$ & $\formatnum{}{0.3513749475869929}$ \\

$(6,3)$ & $2\MultiPartite{\alpha_5n,\, \alpha_5n} \sqcup K_{(1-4\alpha_5) n}$ & $\formatnum{}{0.3650891908417674995}$ & $\formatnum{}{0.3653600282649522}$ \\

$(6,6)$ & 
%$2\MultiPartite{\frac{n}{22}, \frac{n}{22}, \ldots, \frac{n}{22}}$ 
$2\MultiPartite{{n}/{22},\, {n}/{22},\, \ldots,\, {n}/{22}}$
& $\frac{21675}{58564} = \formatnum{}{0.3701079161259476811}$ & $\formatnum{}{0.3701115739691023}$ \\

\hline 

$(7,2)$ & 
%$8K_{\frac{n}{8}}$ 
$8K_{{n}/{8}}$
& $\frac{11025}{32768}=\formatnum{}{0.336456298828125}$ & $\formatnum{}{0.3367351897280793}$ \\

$(7,3)$ & 
%$K_{\frac{3n}{7}} + \O K_{\frac{4n}{7}}$
$K_{{3n}/{7}} \sqcup\O K_{{4n}/{7}}$ 
& $\frac{34560}{117649}=\formatnum{}{0.29375515303997}$ & $\formatnum{}{0.29903663788406776}$ \\

$(7,4)$ & $3\MultiPartite{n/6, \, n/6}$ & $\frac{35}{108} = \formatnum{}{0.324074074074074}$ & $\formatnum{}{0.3269092897814353}$ \\

$(7,5)$ & 
%$3K_{\frac{n}{3}}$ 
$3K_{{n}/{3}}$ 
& $\frac{70}{243}=\formatnum{}{0.2880658436213991}$ & $\formatnum{}{0.2965188292931691}$ \\

$(7,7)$ & 
%$K_{\frac{n}{3}} \sqcup \MultiPartite{\frac{n}{3}, \frac{n}{3}}$ 
$\MultiPartite{\alpha_6n, \, \alpha_6n, \, (1-2\alpha_6)n}$ & $\formatnum{}{0.288086497303598}$ & $\formatnum{}{0.29259270274675236}$ \\

$(7,8)$ & $\MultiPartite{\alpha_7n,\, \alpha_7n,\, \alpha_7n} \sqcup K_{(1-3\alpha_7)n}$ & $\formatnum{}{0.353847617433189}$ & $\formatnum{}{0.35384761752077853}$\\
\hline
\end{tabular}
\caption{The remaining unknown edge-inducibility constants for $\kappa\le 7$. Column 2 describes a (block) construction while Column 3 presents the best lower bound on the edge-inducibility constant coming from this construction (that is, by taking optimal part ratios). Column 4 lists the numerical value returned by SDP solver for $\N=8$.
}
\label{ta:unknown}
\end{center}
\end{table}

There are a few problems where the numerical bound coming from flag algebras is very close to a known lower bound. The most interesting such case is probably the $(5,5)$-edge-inducibility problem, where we conjecture that $\obj[5,5]{}=45/128$, with the lower bound coming from the uniform blowups of the 4-vertex path with loops at end-points. This problem is self-complementary but, unlike the $(4,3)$-case when there are two complementary constructions, here we happen to have only one construction (which is self-complementary). We were not able to round the floating-point matrices returned by the solver into rational matrices that prove the sharp upper bound of $45/128$. One difficulty is that working with $\N=7$ is not enough so one has to use 8-vertex flags (obtaining an SDP program with $|\C F_8^0|=12346$ linear constraints). Also, the rounding step is complicated since there are so-called ``phantom" pairs, namely the pairs witnessing that the problem is not flip-averse (such as the pairs across the two parts for the $(4,3)$-edge-inducibility problem). By looking at the lower order terms, one can show that phantom pairs force some further coefficients $c_F$ to be zero and further vectors to be zero eigenvectors of some matrices $X^\tau$ in~\eqref{eq:FAMain}, when compared to the standard restrictions coming from evaluating  this identity 
on the blowup construction (see~\cite[Equation~(3.10)]{PikhurkoVaughan13} or \cite[Lemma~16]{FalgasMarchantPikhurkoVaughan15}). The presence of such pairs seems to make rounding harder. 
%Another case with very small gap is the inducibility of the graph $F=(5,\{01,02\})$, which is a 3-vertex path plus two isolated vertices. 

The $(7,8)$-edge-inducibility problem is another instance where the bounds are close to each other. Here $\N=7$ does not suffice and, of course, the fact that the optimal part ratios involve the irrational number $\alpha_7$ (a root of an irreducible polynomial of degree 4) makes this task even more challenging.

Further motivation for studying the edge-inducibility problem comes from the question of V.T.S\'os~\cite{Sos12} to describe graphons that are \emph{size forcible}, that is, are determined up to weak isomorphism by the vector of the densities of $(\kappa,\ell)$-configurations over all $(\kappa,\ell)$. We refer the reader to~\cite{Sliacan15,Csoka16,CooleyKangPikhurko22} for the precise definitions and known results.
% Noga Alon (unpublished, see~\cite{Csoka16}) and independently Jakub Sliacan~\cite{Sliacan15} proved that the constant~$\frac12$ graphon is in the Sos family. Then Cs\'oka~\cite{Csoka16} and independently Jacob Fox, Tomasz {\L}uczak and Vera T.\ S\'os (see~\cite{CooleyKangPikhurko22}) proved that constant-$p$ graphon is in the Sos family for any $p\in(0,1)$
Without going into the theory of graphons (for which we refer the reader to the book by Lov\'asz~\cite{Lovasz:lngl}), let us just observe that if Erd\H os-Simonovits stability holds for the $(\kappa,\ell)$-edge-inducibility problem with respect to a unique minimal blowup pattern $B$ and $\V a\in\I S_{m}$ is the unique maximiser of $\obj{\blow{B}{\V a}}$
%with the unique (up to automorphism of $B$) asymptotically optimal vector of $\V a$ of ratios, 
then the limit  of $\V a$-blowups of $B$ is size forcible. This is the case in all our new results listed in Table~\ref{ta:known} except for $(\kappa,\ell)=(4,3)$, so they provide many new examples of graphons that are size forcible (in fact, in the much stronger sense that just one pair $(\kappa,\ell)$ suffices). 
%Likewise, the case of $F$-inducibility problem where we proved perfect stability give a  example of a non-trivial graphon that is forcible just by a single subgraph density, which is a much stronger version of \emph{finite forcibility}.
%, such as the $6$-part construction for the case $(\kappa,\ell)=(4,2)$.

In terms of the inducibility constants of $5$-vertex graphs, as far as we see, there are 5 remaining open cases, namely where $E(F)$ is $\{01,12\}$ (2-edge path plus 2 isolated vertices), $\{01,02,03,12\}$ (triangle with a pendant edge plus 1 isolated vertex), $\{01,12,23,34\}$ (4-edge path), $\{01, 02, 03, 04, 12\}$ (triangle with two pendant edges adjacent to the same vertex) and $\{01, 02, 03, 12, 14\}$ (triangle with two pendant edges adjacent to two different vertices). In the first of these cases, when $F=P_3\sqcup 2K_1$, the numerical bound given by  flag algebra comes very close to the lower bound $\obj[F]{}\ge 174/625$ of Even-Zohar and Linial~\EL{Table~2} given by non-uniform blowups of $3K_2$, namely by $\MultiPartite{n/10,\,n/10}\sqcup 2\MultiPartite{n/5,\,n/5}$ as $n\to\infty$. In the remaining four cases, there seems to be a non-zero gap between the upper and lower bounds; we refer the reader to \EL{} for details.

In an ongoing project  with Jun Gao, Jared Le\'on,  Xizhi Liu and Shumin Sun, we started a systematic study of the $F$-inducibility problem for $6$-vertex graphs~$F$.

Recently, Basit, Granet, Horsley, K\"undgen and Staden~\cite{BasitGranetHorsleyKundgenStaden25x} introduced the following problem for a given blue-red edge-coloured graph $H$ with $\kappa$ vertices (and proved a number of initial results on it). For an (uncoloured) graph $G$, let $\Obj[H]{G}$ be the number of \emph{embeddings} of $H$ into $G$, that is injections $V(H)\to V(G)$ that map red (resp.\ blue) edges of $H$ to edges (resp.\ non-edges) of $G$. Note that we do not stipulate anything about the images of the pairs that are not adjacent in~$H$. 
The \emph{$H$-semi-inducibility problem} asks for $\Obj[H]{n}$, the maximum of $\Obj[H]{G}$ over all $n$-vertex graphs $G$. One can consider the normalised version $\obj[H]{n}:=\Obj[H]{n}/\prod_{i=0}^{\kappa-1}(n-i)$.
%, where $n^{(\kappa)}:=n(n-1)\dots(n-\kappa+1)$ is the $\kappa$-th falling power of $n$. 
Thus $\obj[H]{n}$ is the probability that a random injection $V(H)\to V(G)$ is an embedding of $H$ into~$G$. It is easy to show (see \cite[Proposition~3.1]{BasitGranetHorsleyKundgenStaden25x}) that the limit $\obj[H]{}:=\lim_{n\to\infty}\obj[H]{n}$ exists. 
We call $\obj[H]{}$ the \emph{semi-inducibility constant} of $H$. 
%Note that our counting parameter $\Obj[H]{G}$ is slightly different (up to a scaling factor) from the one used in~\cite{BasitGranetHorsleyKundgenStaden25x} where the authors count the pairs of sets of edges and non-edges in $G$ which form a copy of~$H$. 
\hide{
Of course, each such pair corresponds to the same number of embeddings of $H$ into $G$ (which is the number of colour-preserving automorphisms of $H$), so it is easy to change between these two parameters. We use $\Obj[H]{G}$ as it admits a natural density version $\obj[H]{G}$ and is a more convenient parameter to use in the proof Lemma~\ref{lm:CompletePartite}.
}%
%By symmetry, it is enough to consider edge-coloured $H$ with at least as many blue edges as red edges. 
If $H$ is a colouring of the complete graph on $V(H)$ then the $H$-semi-inducibility problem is, apart from a scaling factor, exactly the inducibility problem for the red subgraph of $H$, which was the main motivation  in~\cite{BasitGranetHorsleyKundgenStaden25x} for introducing this question. 
\hide{
Basit et al~\cite{BasitGranetHorsleyKundgenStaden25x} obtained a number of results on the semi-inducibility problem, and we refer the reader to that paper for further details. 
}

Using flag algebras, we also determined in~\cite{BodnarPikhurko25semi4} the semi-inducibility constant $\obj[H]{}$ for all $H$ with $4$ vertices except when $H$ is the 3-edge path coloured blue-blue-red in this order (or is equivalent to this case). Also, the first author obtained a number of results for 5-vertex graphs~$H$;
%By automating the process, the first author was able (as of now) to resolve the semi-inducibility problem for $5$-vertex $H$ in 193  out of 330 non-equivalent cases. 
these findings will be presented in a follow-up paper. Here we present only the following two results, as they address two open problems highlighted by Basit et al~\cite[Problems~9.1 and 9.2]{BasitGranetHorsleyKundgenStaden25x}.

\begin{theorem}\label{th:Alt6Cycle}
 Let $H$ be the 6-edge cycle with its edges alternatively coloured blue and red as we go along the cycle. Then $\obj[H]{}=1/2^6$.
 \end{theorem}

\bpf
%[Proof of Theorem~\ref{th:Alt6Cycle}]
 The upper bound can be obtained from flag algebras with $\N=6$ (with the certificate named \texttt{semiind\_c6}) while the lower bound comes from $R(K_n,1/2)$, that is, from $1/2$-quasirandom graphs.\epf
 
\begin{theorem}\label{th:Alt3Path}
 Let $H$ be the $3$-edge path with its edges coloured blue, red and blue in this order. Then $\obj[H]{}=2^2/3^3$. Moreover, every sequence of almost extremal graphs is almost $1/3$-regular.
 %Moreover, for every $\e>0$ there are $\delta>0$ and $n_0$ such that if $G$ is a graph with $n\ge n_0$ vertices and $\obj[H]{G}\ge \obj[H]{}-\delta$ then at least $(1-\e)n$ vertices $u$ of $G$ satisfy $\deg_G(u)=(1/3\pm\e)n$.
 \end{theorem}

%Let us now turn to our semi-inducibility results. First, we look at the 3-edge path coloured blue-blue-red.

\bpf
%[Proof of Theorem~\ref{th:Alt3Path}]
The upper bound can be proved via flag algebras with $\N=4$ (with the certificate named \texttt{semiind\_p4}) while the lower bound comes from almost $1/3$-regular graphs. 

Furthermore, our scripts verify that if $\tau$ is any
of the 2-vertex types $\tau_0:=(2,\{\})$ and $\tau_1:=(2,\{01\})$ (the non-edge and the edge), then every vector $\V v$ in the kernel of $X^\tau$ satisfies that
 \begin{equation}\label{eq:Degree}
 \sum_{F\in\C F_3^{\tau}} \left(\p(F_0,F)-\p(F_1,F)\right)\V v_F=0,
 \end{equation}
  where, for $i\in [2]$, $F_i$ denotes the 3-vertex $\tau$-flag with the (unique) unlabelled vertex adjacent to root $i$ but not to root $1-i$. (Thus each coefficient $\p(F_i,F)$ in~\eqref{eq:Degree} is 0 or 1.) 
  
Take any almost  $\obj[H]{}$-extremal graph $G$ of order $n\to\infty$. For  any two distinct vertices $u_0,u_1\in V(G)$, let $\tau$ be $\tau_0$ if $u_0$ and $u_1$ are non-adjacent and let $\tau$ be $\tau_1$ otherwise. If we evaluate the left-hand side of~\eqref{eq:Degree} when  $\V v$ is the vector $\V v_{G,(u_0,u_1)}^{\tau,s}$ of the densities of 3-vertex $\tau$-flags in the $\tau$-flag $G$ rooted at $(u_0,u_1)$, as defined in~\eqref{eq:TauVector}, then we obtain $(\deg_G(u_0)-\deg_G(u_1))/(n-2)$. By  Lemma~\ref{lm:Eigenspaces},  for most choices of $u_0,u_1\in V(G)$, the vector  $\V v_{G,(u_0,u_1)}^{\tau,s}$ is $o(1)$-close to the null-space of $X^\tau$; then the normalised degrees of $u_0$ and $u_1$ must be $o(1)$-close to each other. It follows that $G$ is  almost $\alpha$-regular for some $\alpha=\alpha(n)$ (that potentially depends on~$n$).
  Since $G$ has $(\alpha+o(1))\binom{n}2$ edges and, for most of these edges, both endpoints have the complementary degree $(1-\alpha+o(1))n$, it holds that $\obj[H]{G}=\alpha(1-\alpha)^2+o(1)$. This implies that $\alpha(n)=1/3+o(1)$, that is, extremal graphs are almost $1/3$-regular, as claimed.\epf

Note that every almost $1/3$-regular graph $G$ satisfies $\obj[H]{G}=2^3/3^6+o(1)$, so the above theorem gives a  characterisation of almost extremal graphs.
More can be said about an extremal graph $G$ of order $n\to\infty$. For example, every vertex $u$ of $G$ must have degree $(1/3+o(1))n$. Indeed, since we cannot increase $\Obj[H]{G}$ by replacing $u$ by a clone of any other vertex of $G$, it must be the case that number of embeddings of $H$ that use $u$ is at least $4\Obj[H]{G}/n-12n^2=16n^3/27+O(n^2)$. On the other hand, by the almost $1/3$-regularity of $G$ this number is  equal to 
\[
2\cdot (n-\deg(u))\cdot \frac 13\,n\cdot \frac23\,n  + 2\cdot \deg(u)\cdot (n-\deg(u))\cdot \frac 23\,n+o(n^3),
\]where the first (resp.\ second) term counts those embeddings when $u$ corresponds to one of the two endpoints (resp.\  interior points) of the alternating $3$-edge path. The obtained inequality implies that $d=(1/3+o(1))n$. 
By bootstrapping and analysing more carefully the effect of various local changes on the objective function, further more precise results about the degrees and the structure of any extremal graph should be possible.
%Finally, let us provide formal details about alternating 6-cycle.
%Some partial results can be proved about almost extremal graphs (such as almost $1/2$-regularity). However, we have not been able yet to prove that all almost extremal graphs in Theorem~\ref{th:Alt6Cycle} are $1/2$-quasirandom. 

Since this paper is already rather long, we decided to restrict ourselves here to the above semi-inducibility results.

Independently of this work, the semi-inducibility problem for the alternating path of any odd length was resolved by Chen, Clemen and Noel~\cite{ChenClemenNoel25x} using entropy while the semi-inducibility constant for the alternating 6-cycle was determined by Chen and Noel~\cite{ChenNoel25x}; see also Balogh, Lidicky, Mubayi, Pfender and Volec~\cite{BaloghLidickyMubayiPfenderVolec}.
% whose Theorem 1.1 also imply Theorem~\ref{th:Alt3Path}.
%Also, the semi-inducibility problem  was independently studied using flag algebras by J\'ozsef Balogh, Bernard Lidick\'y, Dhruv Mubayi and Florian Pfender (personal communication).

\section*{Acknowledgements}

The authors thank the anonymous referee for many useful comments.
Both authors were supported by  ERC Advanced Grant 101020255.

\bibliography{bibexport}
%\bibliography{oleg,sets,misc,ramsey,enum,number,posets,sat,ex,matroid,design,random,graph,general,geometry,algorithm,Analysis,limits}

\end{document}